\documentclass[a4paper,11pt, oneside,notitlepage]{article}   
\usepackage{amsmath,amssymb,tocloft,setspace,parskip}       
\usepackage{graphicx}	
\usepackage{float}
\usepackage{amssymb}
\usepackage[natbibapa]{apacite}
\usepackage{booktabs}
\usepackage{rotating}
\usepackage{multicol}
\usepackage[strings]{underscore}
\usepackage{amsthm}
\usepackage{authblk}
\usepackage{enumitem}

\usepackage{authblk}
\title{Principles of Bayesian Inference using General Divergence Criteria}
\author[1,*]{Jack Jewson}
\author[1]{Jim Q Smith}
\author[2]{Chris Holmes}
\affil[1]{Department of Statistics, University of Warwick, Coventry, CV4 7AL}
\affil[2]{University of Oxford, Oxford, OX1 3PA}
\affil[*]{Correspondence address J.E.Jewson@warwick.ac.uk}
\date{February 2018}                     
\begin{document}
\maketitle

\abstract{When it is acknowledged that all candidate parameterised statistical models are misspecified relative to the data generating process, the decision maker (DM) must currently concern themselves with inference for the parameter value minimising the KL-divergence between the model and the process \citep{walker2013bayesian}. However, it has long been known that minimising the KL-divergence places a large weight on correctly capturing the tails of the sample distribution. As a result the DM is required to worry about the robustness of their model to tail misspecifications if they want to conduct principled inference. In this paper we alleviate these concerns for the DM. We advance recent methodological developments in general Bayesian updating \citep{bissiri2016general} to propose a statistically well principled Bayesian updating of beliefs targeting the minimisation of more general divergence criteria. We improve both the motivation and the statistical foundations of existing Bayesian minimum divergence estimation \citep{hooker2014bayesian, ghosh2016robust}, allowing the well principled Bayesian to target predictions from the model that are close to the genuine model in terms of some alternative divergence measure to the KL-divergence. Our principled formulation allows us to consider a broader range of divergences than have previously been considered. In fact we argue defining the divergence measure forms an important, subjective part of any statistical analysis, and aim to provide some decision theoretic rational for this selection.  We illustrate how targeting alternative divergence measures can impact the conclusions of simple inference tasks, and discuss then how our methods might apply to more complicated, high dimensional models.} 




\section{Introduction}\label{Sec:Intro}

In the $M$-closed world, where a fitted model class is implicitly assumed to contain the sample distribution from which the data came, Bayesian updating is highly compelling. However, modern statisticians are increasingly acknowledging that their inference is taking place in the $M$-open world \citep{bernardo2001bayesian}. Within this framework we acknowledge that any class of models we choose is unlikely to capture the actual sampling distribution, and is then at best an approximate description of our own beliefs and the underlying real world process. 
In this situation it is no longer possible to learn about the parameter which generated the data and a statistical divergence measure must be specified between the fitted model and the genuine one in order to define the parameter targeted by the inference \cite{walker2013bayesian}. Remarkably in the $M$-open world standard Bayesian updating can be seen as a method which learns a model by minimising the predictive Kullback-Leibler (KL)-divergence from the model from which the data were sampled \citep{walker2013bayesian, bissiri2016general}. Therefore traditional Bayesian updating turns out to still be a well principled method for belief updating provided the decision maker (DM) concerns themselves with the parameter of the model that is closest to the data as measured by KL-divergence \citep{walker2013bayesian}. Recall the KL-divergence of functioning model $f$ from data generating density $g$ is given by

\begin{equation}
d_{KL}(g,f)=\int\log\left(\frac{g(z)}{f(z)}\right)dG(z).
\label{Equ:KL}
\end{equation}

It is well known that when the KL-divergence is large, which is almost inevitable in high dimensional problems, the KL-divergence minimising models can give a very poor approximation if our main interest is in getting the central part of the posterior model uncertainty well estimated. This is because the KL-divergence gives large consideration to correctly specifying the tails of the generating sample distribution - see Section \ref{SubSub:MovingFromKL}. As a result, when the DM acknowledges that they are in the $M$-open world, they have several options currently available to them:
\begin{enumerate}[leftmargin=*,labelsep=4.9mm]
\item Proceed as though the model class does contain the true sample distribution and conduct {\em{a posteriori}} sensitivity analysis.
\item Modify the model class in order to improve its robustness properties.
\item Abandon the given parametric class and appeal to more data driven techniques.
\end{enumerate}
To this we add:
\begin{enumerate}[resume,leftmargin=*,labelsep=4.9mm]
\item Acknowledge that the model class is only approximate, but is the best available, and seek to infer the model parameters in a way that are most useful to the decision maker. 
\end{enumerate}

Method 1 is how \cite{box1980sampling} recommends approaching parametric model estimation and is the most popular approach amongst statisticians. Although it is acknowledged that the model is only approximate, Bayesian inference is applied as though the statistician believes the model to be correct. The results are then checked to examine how sensitive these are to the approximations made. Authors \cite{berger1994overview} provide a thorough review while \cite{watson2016approximate} consider this in a decision focused manner. 

Method 2 corresponds to the classical robustness approach. Within this approach, one model within the parametric class may be substituted for a model providing heavier tails \citep{berger1994overview}. Alternatively different estimators, for example M-estimators \cite{huber1981robust} \cite{hampel2011robust}, see \cite{greco2008robust} for a Bayesian analogue, are used instead of those justified by the model class. Bayes linear methods \cite{goldstein1999bayes}, requiring many fewer belief statements than a full probability specification, form a special subclass of these techniques.

The third possibility is to abandon any parametric description of the probability space. Examples of this solution include, empirical likelihood methods \citep{owen1991empirical}; a decision focused general Bayesian update \citep{bissiri2016general}; Bayesian non-parametric methods; or to appeal to statistical learning methods such as neural networks or support vector machines.

Though alternatives 2 and 3 can be shown to be very powerful in certain scenarios, it is our opinion that within the context described above that a fourth option holds considerable relative merit. In an applied statistical problem the model(s) represents the DM's best guess of the sample distribution that could be applied. The model provides the only opportunity to input not only structural, but quantitative expert judgements about the domain, something that is often critical to a successful analysis - see \cite{lazer2014parable} for example. The model provides an interpretable and transparent explanation about how different factors might be related to each other. This type of evidence is often essential when advising on decisions to policy makers. Often, simple assumptions play important roles in providing interpretability to the model and in particular preventing it from over fitting to any non-generic features contained in any one data set.

For the above reasons an unambiguous statement of a model, however simple, is in our opinion an essential element for much of applied statistics. In light of this, statistical methodology should also be sufficiently flexible to cope with the fact that the DM's model is inevitably an approximation both of their beliefs and the real world process. Currently, Bayesian statistics sometimes struggles to do this well. We argue below that this is because it implicitly minimises the KL-divergence to the underlying process.

By fitting model parameter in a way that can be non-robust, the DM is having to combine their best guess belief model with something that will be robust to the parameter fitting. The DM must seek the best representation of their beliefs about a process in order to make future predictions. However under traditional Bayesian updating they must also give consideration to how robust these beliefs are. This seems an unfair task to ask of the DM. We therefore propose to decouple what the DM believes about the data generating process, from how the DM wishes the model to be fitted. This results in option 4 above.

This fourth option suggest the DM may actually want to explicitly target more robust divergences, a framework commonly known as Minimum Divergence Estimation (MDE), see \cite{basu2011statistical}. Minimum divergence estimation is of course a well-developed field by frequentists, with Bayesian contributions coming more recently. However when the realistic assumption of being in the $M$-open world is considered the currently proposed Bayesian minimum divergence posteriors fail to fully appreciate the principled justification and motivation required to produce a coherent updating of beliefs. A Bayesian cannot therefore make principled inference using currently proposed methods in the $M$-open setting, except in a way that \cite{miller2015robust} describe as ``tend(ing) to be either limited in scope, computationally prohibitive, or lacking a clear justification''. In order to make principled inference it appears as though the DM must currently concern themselves with the KL-divergence. However in this paper we remove this reliance upon the KL-divergence by providing a justification for Bayesian updating minimising alternative divergences, both theoretically and ideologically. Our updating of beliefs does not produce an approximate or pseudo posterior, but uses general Bayesian updating \citep{bissiri2016general} to produce the coherent posterior beliefs of a decision maker who wishes to produce predictions from a model that provide an explanation of the data that is as good as possible in terms of some pre-specified divergence measure. By doing this the principled statistical practice of fitting model parameters to produce predictions is adhered to, but the parameter fitting is done so acknowledging the $M$-open nature of the problem. 

In Section \ref{Sec:Literature} of this article we review the currently available Bayesian minimum divergence estimation techniques and introduce general Bayesian updating \cite{bissiri2016general}. The contributions of this article are then as follows: in Section \ref{Sec:BayesMDE} we identify the inadequacies in the justification provided by the currently available Bayesian MDE methods and use general Bayesian updating to prove that the Bayesian can still do principled inference on the parameters of the model using alternative, more robust divergences to KL-divergence. This theoretical contribution then allows us to propose a  wider variety of divergences that the Bayesian could wish to minimise than have currently been consider in the literature in Section \ref{Sec:Divergences}. In this section we also consider decision theoretic reasons why targeting alternative divergences to the KL-divergence can be more desirable. Lastly in Section \ref{Sec:Illustrations} we demonstrate the impact model misspecifications can have on a traditional Bayesian analysis for simple inference, regression and time series analysis, and that superior robustness can be obtained by minimising alternative divergences to the KL-divergence. In appendix \ref{Sub:Efficiency} we also show that when the observed data is in fact generated from the model, these methods can be shown not to lose too much precision. In this paper we demonstrate that this advice is not simply based on expedience but has a foundation in a principled inferential methodology. In this purpose we have deliberately restricted ourselves to simple demonstrations designed to provide a transparent illustration of the impact that changing the divergence measure can have on inferential conclusions. However we discuss how robust methodology becomes more important as problems and models become more complex and high dimensional and thus encourage practitioners to experiment with this methodology in practice.

\section{Related work}\label{Sec:Literature}

\subsection{Bayesian theory and scoring functions}

Authors \cite{bernardo2001bayesian} consider Bayesian inference through the lens of a Bayesian decision problem, where quoting a probability belief distribution for future uncertainties is the action to be taken. In this scenario, a scoring rule $S(z,F)$ is used to define the loss of quoting a distribution $F$ and observing a realisation $z$. They argue that the scoring rule associated with scoring probabilistic predictions should be proper and local. A proper scoring rule results in the DM's expected loss being minimised when the DM quotes their true beliefs, and a local scoring rule is one where the score only depends on the quoted probability of the actual observed outcome and nothing else. The only proper, local scoring rule is the logarithmic scoring rule:

\begin{equation}
\ell(f(\cdot;\theta),z)=\sum_{i=1}^n S(z_i,F)=\sum_{i=1}^n-\log\left(f(z_i;\theta)\right).
\end{equation}

The expected additional loss when quoting distribution $F$, with density $f(\cdot;\theta)$, for future outcomes when the data is distributed according to distribution $G$, with density $g$, is then:

\begin{equation}
\mathbb{E}_{Z\sim G}[-\log(f(Z;\theta))]-\mathbb{E}_{Z\sim G}[-\log(g(Z))]=d_{KL}(g,f_\theta),
\end{equation}
the Kullback-Leibler (KL) divergence between the two probability distributions in equation (\ref{Equ:KL}). In fact authors \cite{grunwald2004game} define any divergence $D(\cdot,\cdot)$, associated with proper scoring rule $S(Z,F)$ as the extra penalty for believing $Z$ was distributed according to $F$ when it was actually distributed according to $G$.

\begin{equation}
D(G,F)=S(G,F)-S(G,G)=\mathbb{E}_{Z\sim G}[S(Z,F)]-\mathbb{E}_{Z\sim G}[S(Z,G)].
\end{equation}

\subsection{The data generating process and the \textit{M}-open world}

Sometimes in this paper we use the phrase ``the data generating process''. The data generating process is a widespread term in the literature and appears to suggest that `Nature' is using a simulator to generate observations. While this may fit nicely with some theoretical contributions, it becomes difficult to argue for in reality. In this article we consider the data generating process to represent the DM's true beliefs about the sample distribution of the observations. However in order to correctly specify these the DM must be able to take the time and infinite introspection to consider all of the information available to them/in the world in order to produce probability specifications in the finest of details. As is pointed out by \cite{goldstein1990influence} this requires many more probability specifications to be made at a much higher precision than any DM is ever likely to be able to manage within time constraints of the problem. As a result these genuine beliefs must be approximated. This defines the subjectivist interpretation of the $M$-open world we adopt in this paper - the model used for the belief updating is only ever feasibly an approximation of the DM's true beliefs about the sample distribution they might use if they had enough time to fully reflect. In the special case when the data is the result of a draw from a known probability model - a common initial step in validating a methodology - then this thoughtful analysis and "a data generating process" obviously coincide. Henceforth we use the ``the data generating process'' in this sense to align our terminology as closely as possible with that in common usage.

\subsection{General Bayesian updating}
Following the above motivation, \cite{bissiri2016general} consider only conducting inference about the factors of the world/problem that matter to the DM. They produce what they refer to as a general Bayesian update - a coherent method to produce a posterior distribution over some quantity without relying on a full model for the observations - when considering the non-inferential Bayesian decision problem. Similarly to Bayes linear methods \citep{goldstein1999bayes}, they consider producing posterior beliefs which do not require conditioning. They consider the Bayes act associated with beliefs corresponding to the data generating distribution $G(\cdot)$ and quantity of interest $\theta$ as

\begin{equation}
\theta^{*}=\textrm{arg}\min_{\theta} \mathbb{E}_{Z\sim G}[\ell(\theta,Z)]=\textrm{arg}\min_{\theta}\int \ell(\theta,z)dG(z).
\label{Equ:OptLosDec}
\end{equation}

Rather than elicit a model representing the DM's beliefs over $Z$, \cite{bissiri2016general} argue that given a prior, a loss function and some data $\mathbf{x}$, an updating of beliefs about the value of $\theta^{*}$ must be possible in the absence of a model for the sampling distribution.  They suggest that the posterior distribution resulting from such an updating of beliefs, be as `close' as possible in terms of expected loss to both the prior and the observed data $\mathbf{x}$, used for the updating. In this development the KL-divergence was the only way to measure `closeness' to the prior which resulted in an additive Bayesian update - one where the order the data arrives in does not affect that parameter updating and is consistent with the likelihood principle. In the absence of a likelihood, the best choice to measure the `closeness' of the posterior to the data is to use the expected loss of the observed data under the proposed posterior, \cite{bissiri2016general} show that the posterior minimising the expected loss criteria is:

\begin{equation}
\pi(\theta|\mathbf{x}) \propto \pi(\theta)\exp(-w\ell(\theta,\mathbf{x}))
\label{Equ:GeneralBayesPosterior}
\end{equation}

This posterior is not an approximation. Neither is it a pseudo-posterior, but rather a valid, coherent representation of subjective uncertainty in the minimizer of the decision in equation (\ref{Equ:OptLosDec}). Often $\ell(\theta,x)=\sum_{i=1}^nl(\theta,x_i)$, the cumulative loss of the current data set. Updating using the cumulative loss amounts to replacing integrating over the data generating distribution in equation (\ref{Equ:OptLosDec}), with empirical integration over data whose distribution is $G(z)$. They also mention the need to calibrate the loss-function used in the general Bayesian update against the prior, by setting the value of $w$ in equation (\ref{Equ:GeneralBayesPosterior}). While any probability density is defined so that it integrates to 1, any loss-function can be defined to be arbitrarily large or small. It is important that the weight the loss-function has in the updating process is calibrated against the prior such that the posterior distributions resulting from the general Bayesian update retain probabilistic meaning in the absence of a model. So it follows that the general Bayesian must be prepared to post odds according to their posterior similarly to any other Bayesian. 

It is straightforward to see that if the logarithmic score is used in the general Bayesian update then the traditional Bayesian update is recovered:

\begin{equation}
\pi(\theta|\mathbf{x})\propto \exp(-\sum_{i=1}^n-\log(f(x_i;\theta)))\cdot \pi(\theta)= \prod_{i=1}^nf(x_i;\theta)\pi(\theta).
\end{equation}
This echoes the well-known result that the Bayesian predictive distribution finds the distribution that is closest to the data generating distribution in terms of KL-divergence.
General Bayesian updating allows the goals of the statistical analysis to be coupled together with the parameter updating. This is not something previously possible under traditional Bayesian statistics. While \cite{bissiri2016general} provide the framework to update a probability belief distribution using a loss function, they only consider the logarithmic score in an inferential scenario.

\subsection{M-closed Bayesian minimum divergence estimation}\label{Sub:MclosedBayesMDE}

Minimum divergence estimation (MDE) is an approach to robust inference, largely overlapping with minimum scoring rule inference \cite{dawid2016minimum}. MDE considers making inference about the parameters $\theta$ of a parametric model $\{f(x;\theta):\theta\in\Theta\}$ by minimising the divergence between $f(x;\theta)$ and the data generating function of the observed data. While \cite{dawid2016minimum} uses proper scores, MDE traditionally minimises a member of the class of discrepancies or disparities. Disparity based methods first build a non-parametric density estimate $g_n(x)$ of the data generating distribution from the data (often a Kernel Density estimate). They then conduct inference by minimising the divergence between the model and $g_n(x)$. The frequentist literature in this area is vast so in this article we choose to focus on the Bayesian contributions. These have come first from \cite{hooker2014bayesian} and more recently from \cite{ghosh2016robust} and \cite{ghosh2017general}. 

Authors \cite{hooker2014bayesian} consider posterior densities of the form 

\begin{equation}
\pi(\theta|\mathbf{x})\propto \pi(\theta)\exp(-nd^2_H(g,f(\cdot;\theta))),
\end{equation}

where $g$ is the data generating density and $d_H^2(g,f)$ is the squared Hellinger divergence between $g$ and $f$

\begin{equation}
d_H^2(g,f)=\frac{1}{2}\int(\sqrt{g(z)}-\sqrt{f(z)})^2dz=1-\int\sqrt{g(z)f(z)}dz.
\label{Equ:Hellinger}
\end{equation}

Posteriors of this form are justified following the asymptotic approximation of the divergence to the KL-divergence when the data comes from the model. This results from the fact that the divergence between two probability distributions must be uniquely minimised to 0 when the probability distributions are equal \citep{grunwald2004game}. Therefore, when $g(x)=f(x;\theta_0)$ and as $n\rightarrow\infty$ the distribution minimising any divergences between it and the data generating process will be $f(x;\theta_0)$.

Rather than considering a discrepancy, \cite{ghosh2016robust} considers minimising the proper Tsallis score associated with the density power divergence in order to produce a `pseudo-posterior'. The density power divergence between $g$ and $f$ is 

\begin{equation}
d^{\alpha}_{DPD}(g,f)=\frac{1}{\alpha+1}\int f^{1+\alpha}(z)dz-\frac{1}{\alpha}\int f^{\alpha}(z)g(z)dz +\frac{1}{\alpha(1+\alpha)}\int g^{1+\alpha}(z)dz.
\label{Equ:DensityPowerDivergence}
\end{equation}

This can be seen as a Bayesian analogue of the minimum proper score inference considered by \cite{dawid2016minimum}. Authors \cite{ghosh2016robust} identify the fact that the score associated with the density power divergence does not require a density estimate, but can still be used to provide robust inference. They term the distribution produced by their estimation method a `pseudo-posterior'. However \cite{ghosh2017general} demonstrate exponential convergence results illustrating that this posterior is asymptotically optimal in exactly the same exponential rate as the traditional Bayesian posterior when the model is correctly specified \citep{barron1988exponential}.

\section{Extending Bayesian MDE for the M-open world}\label{Sec:BayesMDE}

The current Bayesian minimum divergence technology is currently only able to produce approximate or pseudo posteriors. Therefore, Bayes rule is currently the only method explicitly available in the literature that DMs can use to produce a well principled updating of their beliefs about the parameters of a model. We address this issue in this section.

\subsection{Why the current justification is not enough}

The equivalence between the KL-divergence and the Hellinger-divergence of the model from the data generating process mentioned in Section \ref{Sub:MclosedBayesMDE} can be seen by taking Taylor expansions about the KL-divergence minimising parameter $\hat{\theta}^{KL}=\arg\min_{\theta\in\Theta}d_{KL}(g(x),f(\cdot;\theta))$. We have that

\begin{align}
\int g_n(x)\log\left(\frac{g_n(x)}{f(x;\theta^{'})}\right)dx&=\int g_n(x)\log\left(\frac{g_n(x)}{f(x;\hat{\theta}^{KL})}\right)dx-(\theta^{'}-\hat{\theta}^{KL})\int\frac{\nabla_{\theta}f(x;\hat{\theta}^{KL})}{f(x;\hat{\theta}^{KL})}g_n(x)dx\nonumber\\
&-\frac{(\theta^{'}-\hat{\theta}^{KL})^2}{2}\int\left(\frac{\nabla_{\theta}^2f(x;\hat{\theta}^{KL})}{f(x;\hat{\theta}^{KL})}-\frac{(\nabla_{\theta}f(x;\hat{\theta}^{KL}))^2}{f(x;\hat{\theta}^{KL})^2}\right)g_n(x)dx+\cdots.\label{Equ:TaylorKL}\\
\int \left(1-\frac{\sqrt{f(x;\theta^{'})}}{\sqrt{g_n(x)}}\right)g_n(x)dx&=\int \left(1-\frac{\sqrt{f(x;\hat{\theta}^{KL})}}{\sqrt{g_n(x)}}\right)g_n(x)dx-(\theta^{'}-\hat{\theta}^{KL})\int\frac{1}{2}\frac{\nabla_{\theta}f(x;\hat{\theta}^{KL})}{\sqrt{f(x;\hat{\theta}^{KL})}\sqrt{g_n(x)}}g_n(x)dx\nonumber\\
-\frac{(\theta^{'}-\hat{\theta}^{KL})^2}{2}\int\frac{1}{2}&\left(\frac{\nabla^2_{\theta}f(x;\hat{\theta}^{KL})}{\sqrt{f(x;\hat{\theta}^{KL})}\sqrt{g_n(x)}}-\frac{1}{2}\frac{(\nabla_{\theta}f(x;\hat{\theta}^{KL}))^2}{\sqrt{g_n(x)}(f(x;\hat{\theta}^{KL}))^{3/2}}\right)g_n(x)dx+\cdots.\label{Equ:TaylorHell}
\end{align}

Now following the same arguments as used in \cite{hooker2014bayesian}, if $g_n(x)$ is consistent for $g=f(\cdot;\hat{\theta}^{KL})=f(\cdot;\theta_0)$ 
then equations (\ref{Equ:TaylorKL}) and (\ref{Equ:TaylorHell}) are equivalent in the limit as $n\rightarrow\infty$. However when $\hat{\theta}^{KL}\neq\theta_0$ is not the data generating parameter because the model class is misspecified, and if $g_n(x)$ is still converging to $g(x)\neq f(x;\theta_0)$ as $n\rightarrow\infty$ then equations (\ref{Equ:TaylorKL}) and (\ref{Equ:TaylorHell}) will be different. In this setting the current literature gives no foundational reasoning why updating using the Hellinger divergence constitutes a principled updating of beliefs.

The `pseudo-posterior' of \cite{ghosh2016robust} is justified as a valid posterior by \cite{ghosh2017general} as the traditional posterior originating from an alternative belief model. However, this alternative model will not even correspond to the DM's approximate beliefs about the data generating process. There is therefore, a lack of formal justification for a DM to update beliefs using Bayes rule on this object.

\subsection{Principled justification for KL in $M$-open world}\label{Sub:PrincipledKL}

In contrast, \citep{walker2013bayesian, bissiri2016general} has provided a principled justification for Bayes' rule in the $M$-open world. This requires the DM to think about the KL-divergence minimising parameter. When the model is correctly specified, the Bayesian learns about the parameter $\theta_0$ that generated the data. This, for any statistical divergence $d(\cdot, \cdot)$ is equivalent to learning the parameter, $\theta_0$, minimising $d(g(\cdot),f(\cdot;\theta))$. When the model is considered to be incorrect, there is no longer any formal relationship between the parameter $\theta$ and the data. The likelihood no longer represents the probability of the observed data conditioned on the parameter. Therefore in order for $M$-open inference to be meaningful a divergence measure must be chosen and the parameter of interest can then defined as $\theta^{*}=\arg\min_{\theta\in\Theta}d(g(\cdot),f(\cdot;\theta))$. Authors \cite{walker2013bayesian} then state that once the parameter of interest has been defined as the minimum of some divergence, it is then possible for a practitioner to define their prior beliefs about where the minimiser of this divergence may lie.  Then the practitioner's final task is to ensure that their Bayesian learning machine is learning about the same parameter with which they defined their prior belief. 

Recalling the well-known result that Bayesian updating learns the parameters of the model which minimises the KL divergence of the model from the data generating density. This allows the DM to continuing conducting belief updating in a principled fashion using Bayes rule provided they are interested in, and specify prior beliefs about, the parameter $\theta^{KL}$  minimising the KL divergence of the model from the data generating process. 

\subsubsection{Moving away from KL-divergence in the M-open world}\label{SubSub:MovingFromKL}

Viewed conversely, \cite{walker2013bayesian} identified that using Bayes rule to update beliefs using a misspecified model results in the DM being implicitly concerned about producing predictions that are closest to the data generating process in terms of KL-divergence. 
In many scenarios using the KL-divergence is the most expedient thing to do. This probably explains why the KL-divergence is so prolific for statistical applications. The logarithmic form of the KL-divergence makes it straightforward to calculate the KL-divergence between many exponential families and the link with Bayes rule allows for conjugate prior distributions to be considered.

However as is well known but often forgotten, someone whose probabilities are elicited using the logarithmic score, or the KL-divergence, will have to beware of approximating the probability of an event by 0. Since if they are wrong, they will incur an infinite loss. As a result, accurate probability specification of tail events will be important to the DM using the logarithmic-score, something which is generally speaking extremely important in pure inference problems \citep{bernardo2001bayesian}.

Often in reality inference is being done in order to produce estimates of expected utilities as part of a larger decision making process. While the tail specification is important for inferential procedures and in situations where the losses are unbounded, for example in gambling and odds setting scenarios, many loss functions connected to real decisions may wish to place less impact on rarely occurring observations in the tails. For entirely reasonable practical reasons these rare events are precisely ones which the DM will find hard to accurately elicit \cite{o2006uncertain}, see \cite{winkler1968evaluation} for a demonstration of this in the context of forecasting the probability of precipitation. Therefore, basing inference on such severe penalties can be unstable. For example a treatment regime is likely to put greater weight on correctly treating the majority of the population than worrying about outlying extreme cases. Authors \cite{watson2016approximate} quite rightly point out that loss functions for real problems are often bounded and any method that simulates parameters via MCMC makes this assumption. As will be demonstrated in Section \ref{Sub:TV}, inferential procedures assuming an unbounded loss function such as the logarithmic score, provide no guarantees about performance in more general decision making scenarios.

Once the consequences of using Bayes rule to solve decision problem in the $M$-open world is understood, we believe that DMs may reasonably desire alternative options for parameter updating, that are as well principled as Bayes rule, but places less importance on tail misspecification. In this new era of ``Big data'' it becomes increasingly likely that the model used for inference is misspecified, especially in the tails of the process - see Section \ref{Sub:HighDim}. We believe many DMs would consider it undesirable for models that approximate the distribution of the majority of the data to be disregarded because they poorly fit a few outlying observations.

\subsection{Principled minimum divergence estimation}\label{Sub:PrincipledMDE}

We next take a foundational approach to theoretically justify an updating of beliefs targeting the parameter minimising any statistical divergence between the model and the data generating density. Consider the general inference problem of wanting to find the parameters $\theta$ of the parametric model $\{f(\cdot;\theta):\theta\in\Theta\}$. The model here can be considered as the DM's best guess at the data generating process. We consider it important to continue to use a model even under the acknowledgement that it is misspecified for the reasons outlined in Section \ref{Sec:Intro}. Following \cite{bernardo2001bayesian} and \cite{walker2013bayesian} we consider fitting these parameters in a decision making scenario, minimising some divergence function
 
\begin{equation}
d(g,f(\cdot;\theta))=\mathbb{E}_{x\sim G}[\ell_d(x,f(\cdot;\theta))]-\mathbb{E}_{x\sim G}[\ell_d(x,g)],
\end{equation}

to the data generating process $g(\cdot)$. The goal is therefore to solve the decision problem

\begin{equation}
\theta^{*}=\arg\min_{\theta} d(g(\cdot),f(\cdot;\theta))=\arg\min_{\theta} \int\ell_d(x,f(\cdot;\theta))dG(x).
\label{Equ:DivMin}
\end{equation}

Here the entropy term in the definition of the divergence is removed from the minimisation because it does not depend on $\theta$. In this scenario \cite{walker2013bayesian} forces the Bayesian into concerns associated with non-robustness by shackling them to minimising the KL-divergence in order to use Bayes' rule. However, equation (\ref{Equ:DivMin}) is identical to the setting consider by \cite{bissiri2016general} in equation (\ref{Equ:OptLosDec}), except that now the loss function depends on the parameter through a model so is a scoring rule. Using equation \ref{Equ:GeneralBayesPosterior}, general Bayesian updating therefore gives us the tools to come up with a coherent updating of beliefs that target the parameter minimising the general divergence $d(\cdot,\cdot)$ as 

\begin{equation}
\pi^{(d)}(\theta|\mathbf{x})\propto \pi^{(d)}(\theta)\exp(-\sum_{i=1}^n\ell_d(x_i,f(\cdot;\theta))),
\label{Equ:GeneralBayesDivergence}
\end{equation}

where we have set the calibration weight $w=1$. This will be discussed further in Section \ref{SubSub:Calibration}. In order to stay consistent with \cite{walker2013bayesian}, we use the notation $\pi^{(d)}$ to indicate that the prior and posterior belief distributions correspond to beliefs about the parameter minimising divergence $d(\cdot,\cdot)$. Applying the general Bayesian update in an inferential scenario like this, provides a compromise between purely loss based general Bayesian inference and traditional Bayesian updating, which is conducted completely independently of the problem specific loss function. By continuing to fit a generative model, the Bayesian minimum divergence posterior somewhat separates the inference from the decision making. The inference is still concerned with estimating model parameters based on how well the model's predictions reflect the data generating process. However, the criteria for closeness between predictions and reality takes decision making into account. As a result, the minimum divergence posteriors can be seen as producing inferences that will be suitable for a broad range of loss functions. The DM is no longer required to exactly define the loss function associated with their decision problem at the inference stage. They need only consider broadly how robust to tail misspecifications they want their inference to be in order to define the target divergence (see Section \ref{Sub:DivergenceComparison}). Considering the realistic $M$-open nature of the model class also justifies making this compromise. The general Bayesian update, assumes absolutely no information about the data generating process, while using Bayes rule traditionally\footnote{Prior to the interpretation provided by  \cite{walker2013bayesian} explained in section \ref{Sub:PrincipledKL}} assumes the data generating density is known precisely. As we point out, it is actually more likely that the decision maker is able to express informative but not exact beliefs about the data generating process and therefore a half-way-house between these two is appropriate in reality.

Plugging in the corresponding loss functions gives general Bayesian posteriors minimising the Hellinger and density power divergence as
\begin{align}
\pi^H(\theta|\mathbf{x})&\propto \pi^H(\theta)\exp\left(\sum_{i=1}^n\frac{\sqrt{f(x_i;\theta)}}{\sqrt{g_n(x_i)}}\right)\label{Equ:GeneralBayesH}\\
\pi^{DPD}_{\alpha}(\theta|\mathbf{x})&\propto \pi^{DPD}_{\alpha}(\theta)\exp\left(\sum_{i=1}^n\left\lbrace \frac{1}{\alpha}f(x_i;\theta)^{\alpha}-\frac{1}{1+\alpha}\int f(y;\theta)^{\alpha+1}dy\right\rbrace  \right).
\label{Equ:GeneralBayesDPD}
\end{align}

Equation (\ref{Equ:GeneralBayesH}) introduces $g_n(\cdot)$  to estimate the data generating density $g$ (see Section \ref{Sub:DensityEstimation} for more on this). As a result, the general Bayesian updating is being conducted using an empirical loss function, say $\hat{\ell}_{H}$, which approximates the true loss function associated with minimising the Hellinger divergence between the model and the data generating process $\ell_H$.

Equation (\ref{Equ:GeneralBayesDPD}) is exactly the distribution resulting from the robust parameter update of \cite{ghosh2016robust}, while equation (\ref{Equ:GeneralBayesH}) is similar to the posterior produced by \cite{hooker2014bayesian} except the divergence function has been decomposed into its score and entropy term here. This demonstrates that the posteriors above are not pseudo posteriors - as~\cite{ghosh2016robust} suggests - or approximations of posteriors - as \cite{hooker2014bayesian} suggests - they are the correct way for the DM to update their prior beliefs about the parameter minimising an alternative divergence to the KL-divergence. 

Authors \cite{basu2011statistical} showed that $\hat{\theta}_n=\arg\max_{\theta}\sum_{i=1}^n\left\lbrace \frac{1}{\alpha}f(x_i;\theta)^{\alpha}-\frac{1}{1+\alpha}\int f(y;\theta)^{\alpha+1}dy\right\rbrace $ is a consistent estimate of $\theta^{DPD}$, while \cite{ghosh2016robust} use this to show that under certain regularity conditions $\pi_{\alpha}^{DPD}$ follows a Bernstein-von Mises (BVM) type result, but with convergence in probability, where the asymptotic variance of the posterior is the expected second derivative of the loss function $\ell^{DPD}(x_i,f(\cdot,\hat{\theta}_n))$. Authors \cite{hooker2014bayesian} prove the corresponding asymptotic normality result for the posterior minimising the Hellinger divergence, but assume that the data generating process is contained within the model class, as this was required to justify their posterior in the first place. Consistency to the minimiser of any divergence provides consistency with the data generating process when it is contained within the model class, as by definition all divergences are minimised to 0 when the distributions are the same. The results of \cite{ghosh2017general} demonstrate that not only is the posterior mean consistent, but the whole posterior is asymptotically optimal in the same sense as the traditional Bayesian posterior is when the model class contains the data generating process.

Another principled alternative to traditional Bayesian updating when it is difficult to fully specify a model for the data generating process is Bayes linear methods \citep{goldstein1999bayes}. These only require the subjective specifications of expectations and covariances for various quantities the DM is well informed about and interested in order to do inference. Alternatively by choosing a more robust divergence to use for the updating our method tackles the same problem by being robust to routine assumptions making a full probability specification a much less strenuous task.

\subsubsection{The likelihood principle and Bayesian additivity}\label{Sec:LikelihoodPrinciple}

Authors \cite{hooker2014bayesian} identify that their posterior distribution minimising the Hellinger divergence no longer satisfies Bayesian additivity. That is to say that the posteriors would look different if the data were observed in one go or in two halves for example. This is because the posterior now depends on the estimate of the data generating density, which depends in a non-additive way on the data. This constitutes a departure from the likelihood principle underpinning traditional Bayesian statistics. The likelihood principle says that the likelihood is sufficient for the data. When the model is perfectly specified this is a sensible principle, the likelihood of the observed data under the correct model represents all the information in the data. However in the $M$-open world the likelihood principle is no longer a reasonable requirement. When the model is only considered as an approximation of the data generating process, it is not unreasonable to suspect that the data contains more information than is represented in the likelihood of an incorrect model. However, Bayesian additivity (defined as coherence by \cite{bissiri2016general}) is a central principle of the general Bayesian update. \cite{bissiri2016general} chose the KL-divergence to measure the distance between the prior and the posterior as it was the only divergence measure that left the belief update additive. The general Bayesian posterior minimising the Hellinger divergence combines the loss function for each observation in an additive way, which is consistent with the additivity demanded by \cite{bissiri2016general}, but a different density estimate is produced when the data is considered as a whole or in parts. This causes the empirical loss function used for the updating to be different when the data arrives in sections, as opposed to one go. If the data generating density were available then the exact loss function associated with minimising the Hellinger divergence could be calculated and the Bayesian update would be additive. However because an approximation of the loss function associated with the Hellinger divergence is used, additivity is sacrificed. 

While the posterior of \cite{ghosh2016robust} is still additive, the likelihood is also not considered to be sufficient for the data, as the posterior also depends on the integral, $\int f(y;\theta)^{1+\alpha}dy$. Using an alternative divergence to the KL-divergence in order to conduct Bayesian updating, requires that additional information to the `local' information provided by the likelihood of the observed data is incorporated into the loss function. For the Hellinger-divergence this information comes from the data in the way of a density estimate, while for the density power-divergence this comes from the model through $\int f(y;\theta)^{1+\alpha}dy$.

\subsubsection{A note on calibration}\label{SubSub:Calibration} The two posteriors in equations (\ref{Equ:GeneralBayesH}) and (\ref{Equ:GeneralBayesDPD}) have set the general Bayesian calibration weight to 1. Unlike probability distributions, loss functions can be of arbitrary size and therefore it is important that they are calibrated with the prior. However, this is not a problem here. The posteriors above both use well-defined models and the canonical form of well-defined divergence functions and as a result there is no arbitrariness in the size of the loss. This is further demonstrated by the fact that Bayes rule corresponds to using the canonical form of the KL-divergence and a probability model with weight $w=1$. We do note that when the model is correct the posterior variance of these method will be comparatively bigger for finite data samples than that of the traditional Bayesian posterior. This is to be expected, \cite{zellner1988optimal} showed that Bayes rule processes information optimally and therefore produces the most precise posterior distributions. The posterior is simply a subjective reflection of the DM's uncertainty after seeing the data, and if they believe their model is incorrect and therefore target a more robust divergence, they are likely to have greater posterior uncertainty than if they naively believe their model to be correct and proceed accordingly.

\section{Possible divergences to consider}\label{Sec:Divergences}

The current formulation of Bayesian MDE methods has limited which divergences have been proposed for implementation. However demonstrating that principled inference can be made using alternative divergence measures than the KL-divergence, as we have in Section \ref{Sub:PrincipledMDE}, allows us to consider the selection of the divergence used for updating to be a subjective judgement made by the DM, alongside the prior and model, to help tailor the inference to the specific problem. Authors \cite{celeux2017some} observed that while \cite{gelman2015beyond} advocate greater freedom for subjective judgements to impact statistical methodology, they fail to consider the possibility of subjective Bayesian parameter updating. In the $M$-closed world Bayes rule is objectively the correct thing to do, but in the $M$-open world this is no longer so clear. Very few problems seek answers that are connected with a specific dataset or model, they seek answers about the real world process underpinning these. Authors \cite{goldstein2006subjective}, focusing on belief statements, demonstrate that subjective judgements help to generalise conclusions from the model and the data to the real world process. Carefully selecting an appropriate divergence measure can further help a statistical analysis to do this. We see this as part of the meta-modelling process where the DM can embed beliefs about their nature and extent of their modelling class explicitly. 

Below we list some divergences that a practitioner may wish to consider. We do not claim that the list is exhaustive, but merely contains the divergences we have considered using. We qualitatively describe the features of the inferences produced by targeting each specific divergence, which will be demonstrated empirically in Section \ref{Sec:Illustrations}. Throughout we consider the Lebesgue measure to be the dominating measure.

\subsection{Total Variation Divergence}\label{Sub:TV}

The Total-Variation divergence between probability distributions $g$ and $f$ is given by 

\begin{equation}
d_{TV}(g,f)=\sup_A|g(A)-f(A)|=\frac{1}{2}\int|g(z)-f(z)|dz,
\end{equation} 

and the general Bayesian update targeting the minimisation of the TV-divergence can be produced by using loss function 

\begin{equation}
\ell_{TV}\left(x,f(\cdot;\theta)\right)=\frac{1}{2}\left|1-\frac{f(\cdot;\theta)}{g_n(x)}\right|
\label{Equ:TVloss}
\end{equation}

in equation (\ref{Equ:GeneralBayesDivergence}). It is worth noting however that the loss function in equation (\ref{Equ:TVloss}) is not monotonic in the predicted probability of each observation. This is discussed further in Appendix \ref{Sub:Efficiency}. In fact when informed decision making is the goal of the statistical analysis, closeness in terms of TV-divergence ought to be the canonical criteria the DM demands.  If $d_{TV}(g,f)\leq\epsilon$ then for any utility function bounded by 1, the expected utility under $g$ of making the optimal decision believing the data was distributed according to $f$ is at most $2\epsilon$ worse than the expected utility gained from the optimal decision under $g$ (see \cite{smith2010bayesian} for example). Therefore if the predictive distribution is close to the data generating process in terms of TV-divergence then the consequences of the model misspecification in terms of the expected utility of the decisions made is small. However explicit expressions of the TV-divergence even between known families, are rarely available. This somewhat hinders an algebraic analysis of the divergence. 

It is straightforward to see that $0\leq d_{TV}(g,f)\leq 1$ and the KL-divergence form an upper bound on the TV-divergence   through Pinsker's inequality.
However, the KL-divergence does not bound the TV-divergence below so there are situations where the TV-divergence is very small but the KL-divergence is very large. In this scenario, a predictive distribution whose associated optimal decisions achieve a close to optimal expected utility estimate (as the distribution is close in TV-divergence) will receive very little posterior mass. This is clearly undesirable in a decision making context.

Authors \cite{hooker2014bayesian} identify some drawbacks of having a bounded score function. The score function being upper bounded means that there is some limit to the score that can be incurred in the tails of the posterior distribution. The score incurred for one value of $\theta$, will here be very similar to the score incurred by another value of $\theta$ past a given size. Therefore the tails of the posterior will be equivalent to the tails of the prior. As a result the DM is required to think more carefully about their prior distribution. Not only are improper priors prohibited, but more data is required to move away from a poorly specified prior. This can result in poor finite sample efficiency when the data generating process is within the chosen model class. 



\subsection{Hellinger Divergence}

Authors \cite{devroye1985nonparametric} observed, that the Hellinger divergence can be used to bound the TV-divergence both above and below: 

\begin{equation}
d^2_H(g,f)\leq d_{TV}(g,f)\leq \sqrt{d_H^2(g,f)}\sqrt{2-d_H^2(g,f)}\leq \sqrt{2d^2_H(g,f)}.
\end{equation} 

where 
\begin{equation}
d_H^2(g,f)=\frac{1}{2}\int(\sqrt{g(z)}-\sqrt{f(z)})^2dz=1-\int\sqrt{g(z)f(z)}dz
\end{equation}
As a result, the Hellinger-divergence and the TV-divergence are geometrically equivalent. Thus, if one of them is small, the other is small and similarly if one of them is large the other also is. So if one distribution is close to another in terms of TV-divergence, then the two distributions will be close in terms of Hellinger-divergence as well. Authors \cite{beran1977minimum} first noted that minimising the Hellinger-divergence gave a robust alternative to minimising the KL-divergence, while \cite{hooker2014bayesian} proposed a Bayesian alternative (equation (\ref{Equ:GeneralBayesH})) that has been discussed at length above. While \cite{hooker2014bayesian} motivated their posterior through asymptotic approximations, identifying the geometric equivalence between the Hellinger-divergence and TV-divergence proposes further justification for a robust Bayesian updating of beliefs similar to that of \cite{hooker2014bayesian}. Specifically, if being close in terms of TV-divergence is the ultimate robust goal, then being close in Hellinger-divergence will guarantee closeness in TV-divergence. Hellinger can therefore serve as a proxy for TV-divergence that retains some desirable properties of the KL-divergence: it is possible to compute the Hellinger divergence between many known families \citep{smith1995bayesian} and the score associated with the Hellinger divergence has a similar strictly convex shape, this is discussed further in Appendix \ref{Sub:Efficiency}. 
A posterior targeting the Hellinger-divergence does suffer from the same drawbacks associated with having a bounded scoring function that are mentioned at the end of the previous section. Lastly along with TV-divergence, the Hellinger-divergence is also a metric.

\subsection{$\alpha\beta$-divergences}\label{Sub:alphabetaDivergence}

A much wider class of divergences are provided by the two parameter $\alpha\beta$-divergence family proposed in \cite{cichocki2011generalized} 

\begin{equation}
D_{AB}^{(\alpha,\beta)}(g,f)=\int\frac{1}{\alpha(\alpha+\beta)}f^{\alpha+\beta}(z)-\frac{1}{\alpha\beta} g^{\alpha}(z)f^{\beta}(z)+\frac{1}{\beta(\alpha+\beta)}g^{\alpha+\beta}(z)dz,\quad \alpha,\beta\geq 0,\alpha+\beta\neq0.
\end{equation}

Reparametrising this such that $\alpha=(1+\lambda(1-\alpha_S))$ and $\beta=(\alpha_S-\lambda(1-\alpha_S))$ gives the S-divergence of \cite{ghosh2017generalized}.



The general Bayesian posterior targeting the parameter minimising the $\alpha\beta$-divergence of the model from the data generating density is

\begin{equation}
\pi^{\alpha\beta}(\theta|\mathbf{x})\propto\pi^{\alpha\beta}(\theta)\exp(\sum_{i=1}^n\frac{1}{\alpha\beta} g_n^{\alpha-1}(x_i)f^{\beta}(x_i;\theta)-\int\frac{1}{\alpha(\alpha+\beta)}f^{\alpha+\beta}(z;\theta)dz)
\end{equation}

Letting $g$ represent the data generating density, and $f(\cdot;\theta)$ some predictions of $g$ parametrised by $\theta$, \cite{cichocki2011generalized} provide intuition about how the parameters $\alpha$ and $\beta$ impact the influence each observation $x$ has on the inference about $\theta$. They consider influence to be how the observations impact the estimating equation of $\theta$, this is a frequentists setting but the intuition is equally valid in the Bayesian setting.

\begin{equation}
\begin{cases}
\alpha>1 \textrm{, down weights x with smaller ratios }g(x)/f(x;\theta)\textrm{ with respect to larger ones.} \\
\alpha<1 \textrm{, down weights x with the larger ratios }g(x)/f(x;\theta)\textrm{ with respect to smaller ones.} \\
\beta+\alpha>1 \textrm{, down weights x where } f(x;\theta)\textrm{ is small.} \\
\beta+\alpha<1 \textrm{, down weights x  where } f(x;\theta)\textrm{ is large.} \\
\end{cases}
\end{equation}

The size of $g(x)/f(x;\theta)$ for an observation $x$ defines how outlying (large values) or  inlying (small values) the observation is. Choosing $\alpha<1$ ensures outliers relative to the model have little influence on the inference, while adjusting the value of $\beta$ such that $\beta+\alpha>1$ down-weights the influence of unlikely observations from the model $f(x;\theta)$. 

The hyperparameters $\alpha$ and $\beta$ control the trade-off between robustness and efficiency and selecting these hyperparameters is part of the subjective judgement associated with selecting that divergence. This being the case we feel that to ask for values of $\alpha$ and $\beta$ to be specified by a DM is perhaps over ambitious, even given the interpretation ascribed to these parameters above. However, using the $\alpha\beta$-divergence for inference can provide greater flexibility to the DM, which may in some cases be useful. Here to simplify the subjective judgements made we focus on two important, one parameter subfamilies within the $\alpha\beta$-divergence family - the alpha-divergence where $\beta=1-\alpha$ and the beta-divergence where $\alpha=1$.

\subsubsection{alpha Divergence}

The alpha-divergence, introduced by \cite{csisz1967information} and extensively studied by \citep{shun2012differential} is 
\begin{align}
\begin{split}
d_{\alpha}(g,f)&=\frac{1}{\alpha(1-\alpha)}\left\lbrace 1-\int g(z)^{\alpha}f(z)^{1-\alpha}dz\right\rbrace,
\end{split}
\label{Equ:alphaDiv}
\end{align}
where $\alpha\in\mathbb{R}\backslash\{0,1\}$. There exists various reparametrisations of this: Amari notation uses $\alpha_A$ with $\alpha=\frac{1-\alpha_A}{2}$; or Cressie-Read notation \citep{cressie1984multinomial} introduces $\lambda$ with $\alpha=\lambda+1$.  

We generally restrict attention to values of $\alpha\in(0.5,1)$. $\alpha=1$ corresponds to a KL-divergence limiting case while $\alpha=0.5$ is 4 times the Hellinger divergence.  We think of these as two extremes of efficiency and robustness within the alpha-divergence family with which a DM would want to conduct inference between. The parameter $\alpha$ thus controls this trade-off.

The general Bayesian posterior targeting the minimisation of the alpha-divergence is 

\begin{equation}
\pi^{\alpha}(\theta|\mathbf{x})\propto \pi^{\alpha}(\theta)\exp(\frac{1}{\alpha(1-\alpha)}\sum_{i=1}^ng_n(x_i)^{\alpha-1}f(x_i;\theta)^{1-\alpha}).
\label{Equ:alphaGB}
\end{equation}

Note the power on $g$ in equation (\ref{Equ:alphaGB}) is reduced by 1 from equation (\ref{Equ:alphaDiv}) as we consider empirical expectations in order to estimate the expected score associated with the alpha-divergence. 

It was demonstrated in \cite{sason2015bounds} (Corollary 1) that for $\alpha\in(0,1)$ the alpha-divergence can be bounded above by TV-divergence :
\begin{equation}
\alpha(1-\alpha)d_{\alpha}(g,f)\leq d_{TV}(g,f).
\end{equation}
Therefore, if the TV-divergence is small then the alpha-divergence will be small (provided $\alpha\neq\{0,1\}$). So a predictive distribution that is close to the data generating density in terms of TV-divergence will receive high posterior mass under an update targeting the alpha-divergence. Once again, the general Bayesian posterior requires a density estimate to compute the empirical estimate of the score associated with the divergence and the same issues associated with having a bounded score function that were discussed in Section \ref{Sub:TV} apply here.

Considering the alpha divergence as a subfamily of the $\alpha\beta$-divergence takes $\beta=1-\alpha$, and therefore $\alpha+\beta=1$. Therefore the influence of observations $x$ are weighted based on their ratio of $g(x)/f(x;\theta)$ only and not on the value of $f(x;\theta)$ in isolation. $\alpha=1$ corresponded to the KL-divergence and in order to obtained greater robustness $\alpha$ is chosen such that $\alpha<1$ in order to down-weight the influence larger ratios of $g(x)/f(x;\theta)$, corresponding to outlying values of $x$, have on the analysis. However down-weighting some of these ratios relative to the KL-divergence sacrifices efficiency. This demonstrates the efficiency and robustness trade-off associated with the alpha-divergence.

\subsubsection{beta Divergence}

The beta-divergence, also known as the density power-divergence \cite{basu1998robust}, is a member of the Bregman divergence family
\begin{equation}
D(g,f)=\int \psi\{g(z)\}-\psi\{f(z)\}-\psi{'}\{f(z)\}(g(z)-f(z)) d\mu(z).
\end{equation}

Taking $\psi=t^{\beta+1}$, to be the Tsallis score returns the beta-divergence

\begin{equation}
d^{\beta}(g,f)=\frac{1}{\beta+1}\int f^{\beta+1}(z)dz-\frac{1}{\beta}\int f^{\beta}(z)g(z)dz+\frac{1}{\beta(\beta+1)}\int g^{\beta+1}(z)dz,
\end{equation}
$\beta\in\mathbb{R}\backslash\{-1,0\}$. This results in the density power divergence given in equation (\ref{Equ:DensityPowerDivergence}), parameterised by $\beta$ rather than $\alpha$ here. Both \cite{basu1998robust} and \cite{dawid2016minimum} noticed that inference can be made using the beta-divergence without requiring a density estimate. This was used in \cite{ghosh2016robust} to produce a robust posterior distribution that did not require an estimate of the data generating density, which has been extensively discussed in previous sections. 

Considering the density power divergence as part of the $\alpha\beta$-divergence family takes $\alpha=1$, which results in treating all values of $g(x)/f(x;\theta)$ equally in isolation, but taking $\beta>0$ results in $\alpha+\beta>1$ which means observations that have low predicted probability $f(x;\theta)$ under the model are down-weighted. Under the KL-divergence at $\beta=0$, the influence of an observation $x$ is inversely related to its probability under the model, through the logarithmic score. Raising $\beta$ above 0 will decrease the influence of the smaller values of $f(x;\theta)$, robustifying the inference to tail specification. However this results in a decrease in efficiency relative to methods minimising the KL-divergence \citep{ghosh2017general}. As a result there may be some issues using the density power-divergence when high dimensional observations and models are considered. As the dimension increases the predicted probability of each (multivariate) observation shrinks towards 0. The density power divergence down-weights the influence of observations with small predicted probabilities, and as a result $\beta$ ($\alpha$ when we consider the density power divergence) needs to be selected very carefully in order to prevent the analysis from disregarding the majority of the data. This is a price that is paid in order to not require a data generating density estimate.

\subsubsection{The S-Hellinger divergence}

One further one-parameter special case of the $\alpha\beta$-divergence (S-divergence) is the S-Hellinger divergence given by \cite{ghosh2017generalized}

\begin{equation}
d_{SH}(g,f)=\frac{2}{1+\alpha_S}\int \left(g(z)^{(1+\alpha_S)/2}-f(z)^{(1+\alpha_S)/2}\right)^2dz.
\end{equation}

This is generated from the S-divergence by taking $\lambda=-\frac{1}{2}$. Taking $\alpha_S=0$ recovers twice the squared Hellinger divergence and $\alpha=1$ gives the $L_2$ squared divergence. \cite{ghosh2017generalized} observe that the S-Hellinger divergence is a proper distance metric. Translating these back into the notation of $\alpha\beta$-divergences gives 
\begin{equation}
\alpha=\frac{1}{2}\left(1+\alpha_S\right), \quad \beta=\frac{1}{2}\left(\alpha_S+1\right).
\end{equation}
Therefore $\alpha_s\in[0,1]$ gives $\alpha<1$ and $\alpha+\beta>1$. As a result, the squared Hellinger divergence down-weights the influence of large ratios of $g(x)/f(x;\theta)$ with respect to smaller ones and also down-weights these ratios when $f(x;\theta)$ is small. In fact, there is a trade-off here: $\alpha_s$ moving closer to 1 increase $\alpha+\beta$ away from 1, thus increasing the amount that ratios of $g(x)/f(x;\theta)$ are down-weighted for small $f(x;\theta)$. However as $\alpha_s$ moves closer to 1, $\alpha$ also draws closer to 1, which reduced the amount ratios of $g(x)/f(x;\theta)$ are down weighted with respect to smaller ones. Once again we consider trading these two off against each other to be too greater task for the DM to consider. We therefore do not consider the S-Hellinger divergence again here and instead for brevity consider only the choice between the one parameter alpha and density power divergences. In any case, bridging the gap between the Hellinger and the KL-divergence, as the alpha-divergence does, appears to make more sense in an inferential setting than bridging the gap between the Hellinger and $L_2$ squared divergence.

\subsection{Comparison}\label{Sub:DivergenceComparison}

The bullet points below summaries the reasons, introduced above, for which a DM might prefer a specific divergence. Once again we stress that we consider the choice of divergence to be a subjective and context specific judgement to be made by the DM similarly to their prior and model.

\begin{itemize}[leftmargin=*,labelsep=5.8mm]
\item \textbf{KL-divergence}: suitable if the tail behaviour of the model is of high importance or there is the belief that the tails of the model are exactly correctly specified. Minimising the KL-divergence also makes best use out of the data so may be suitable if the sample size is very small and the prior is not necessarily informative.
\item \textbf{TV-divergence/Hellinger-divergence}: suitable if the data set has many observations and robustness for a decision problem with a bounded loss function is of high importance. The Hellinger-divergence appears to be computationally easier to work with (see Appendix \ref{Sub:Efficiency}). Both require accurate density estimation
\item \textbf{alpha-divergence}: suitable for when a trade-off between robustness and efficiency is considered. Requires access to a density estimation technique. 
\item \textbf{The density power divergence}: suitable again for when a trade-off between robustness and efficiency is considered. Does not require a density estimate but does require $\alpha$ to be set carefully so that the data is not completely disregarded by the updating process.
\end{itemize}

Unfortunately there is no free lunch when it comes to applying any of these methods to complex, high dimensions problems. Minimising the KL-divergence is the most computationally efficient. It provides the possibility of implementing conjugate prior distributions and avoiding the computational burden of producing a density estimate. However as the dimensions and complexity of the problem increases, robustness to model misspecification will become more and more important. Minimising the Hellinger-divergence, TV-divergence and alpha-divergence requires an estimate of the data generating density which is certainly not straightforward for high dimensional data (see the discussion in Section \ref{Sub:DensityEstimation}). However once this has been estimated then $\alpha$ can just be selected based on how important tail misspecifications are, with a guarantee on some reasonable efficiency. 
Minimising the density power divergence has the computational advantage of not requiring an estimate of the data generating density. However when the dimension of the problem is high, much more thought needs to be put into specifying the parameter $\alpha$ to ensure that the updating does not completely disregard the data. The technology for doing this effectively is still in its infancy. Without the data generating density, the density power divergence requires information on the shape of the data from another form in order to fit efficiently and robustly.


\subsection{Density estimation}\label{Sub:DensityEstimation}

As has been mentioned before, for the TV, Hellinger and $\alpha\beta$ divergence, it is not possible to exactly calculate the loss function associated with any value of $\theta$ and $x$ because the data generating density $g(x)$ will not be available. In this case, a density estimate of $g(x)$ is required to produce an empirical loss function. The Bayesian can consider the density estimate as providing additional information to the likelihood from the data (see Section \ref{Sec:LikelihoodPrinciple}'s discussion on the likelihood principle), and can thus consider their general Bayesian posterior inferences to be made conditional upon the density estimate as well as the data. The general Bayesian update is a valid update for any loss function, and therefore conditioning on the density estimate as well as the data still provides a valid posterior. However, how well this empirical loss function approximates the exact loss function associated with each divergence ought to be of interest. The exact loss function is of course the loss function the DM would prefer to use having made the subjective judgement to minimise that divergence. If the density estimate is consistent to the data generating process, then provided the sample size is large the density estimate will converge to the data generating density, and the empirical loss function will then correctly approximate the loss function associated with that divergence. It is this fact that ensures the consistency of the posterior estimates of the minimum Hellinger posterior \cite{hooker2014bayesian}. 

Authors \cite{hooker2014bayesian} use a fixed width kernel density estimate (FKDE) to estimate the underlying data generating density and in our examples in Section \ref{Sec:Illustrations} we adopt this practice using a Radial Basis Function (RBF) kernel for simplicity and convenience. However we note that \cite{silverman1986density} identifies practical drawbacks of FKDEs, including their inability to correctly capture the tails of the data generating process, whilst not over smoothing the centre, as well as the number of data points required to fit these accurately in high dimensions. In addition to this \cite{tamura1986minimum} observe that the variance of the FKDE when using a density kernel in high dimensions lead to asymptotic bias in the estimate that is larger than $\mathcal{O}\left(n^{-1/2}\right)$. Alternatives include using a kernel with better mean-squared error properties \citep{epanechnikov1969non,rosenblatt1976maximal}, variable width adaptive KDEs \citep{abramson1982bandwidth}, which \cite{hwang1994nonparametric} show to  be promising in high dimensions, piecewise-constant (alternatively tree based) density estimation \citep{ram2011density, lu2013multivariate} which are also promising in high dimensions, or a fully Bayesian Gaussian process as is recommended in \cite{li2016framework}. 

\section{Illustrations}\label{Sec:Illustrations}

In this section we aim to illustrate some of the qualitative features associated with conducting inference targeting the minimisation of the different divergences identified in Section \ref{Sec:Divergences}. Throughout these experiments \textit{stan} \citep{carpenter2016stan} is used to produce fast and efficient samples from the general Bayesian posteriors of interest
. 

\subsection{$M$-open robustness}\label{Sub:Robustness}

\subsubsection{Simple inference}

\begin{figure}[H]
\centering
\includegraphics[width=5cm]{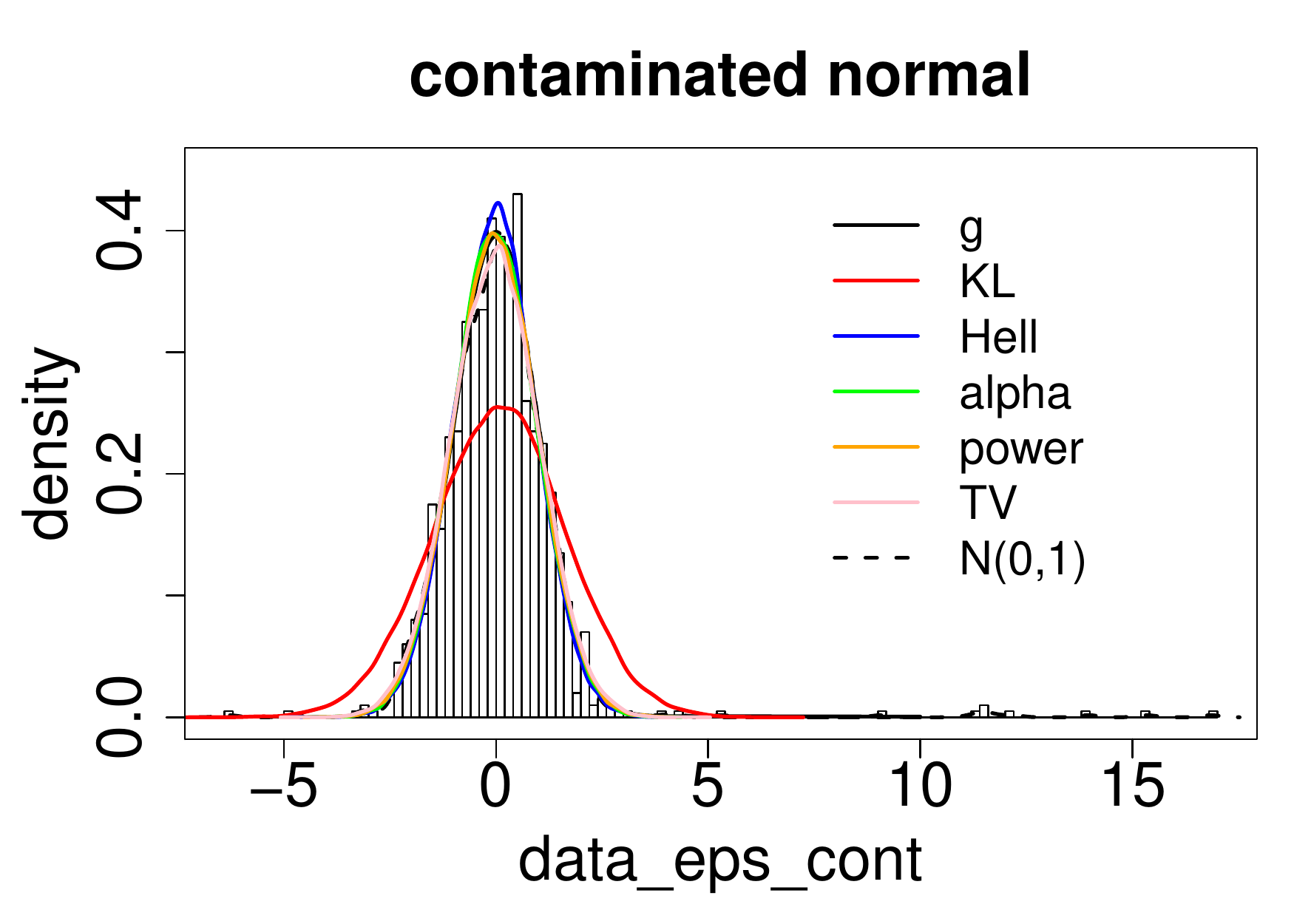}
\includegraphics[width=5cm]{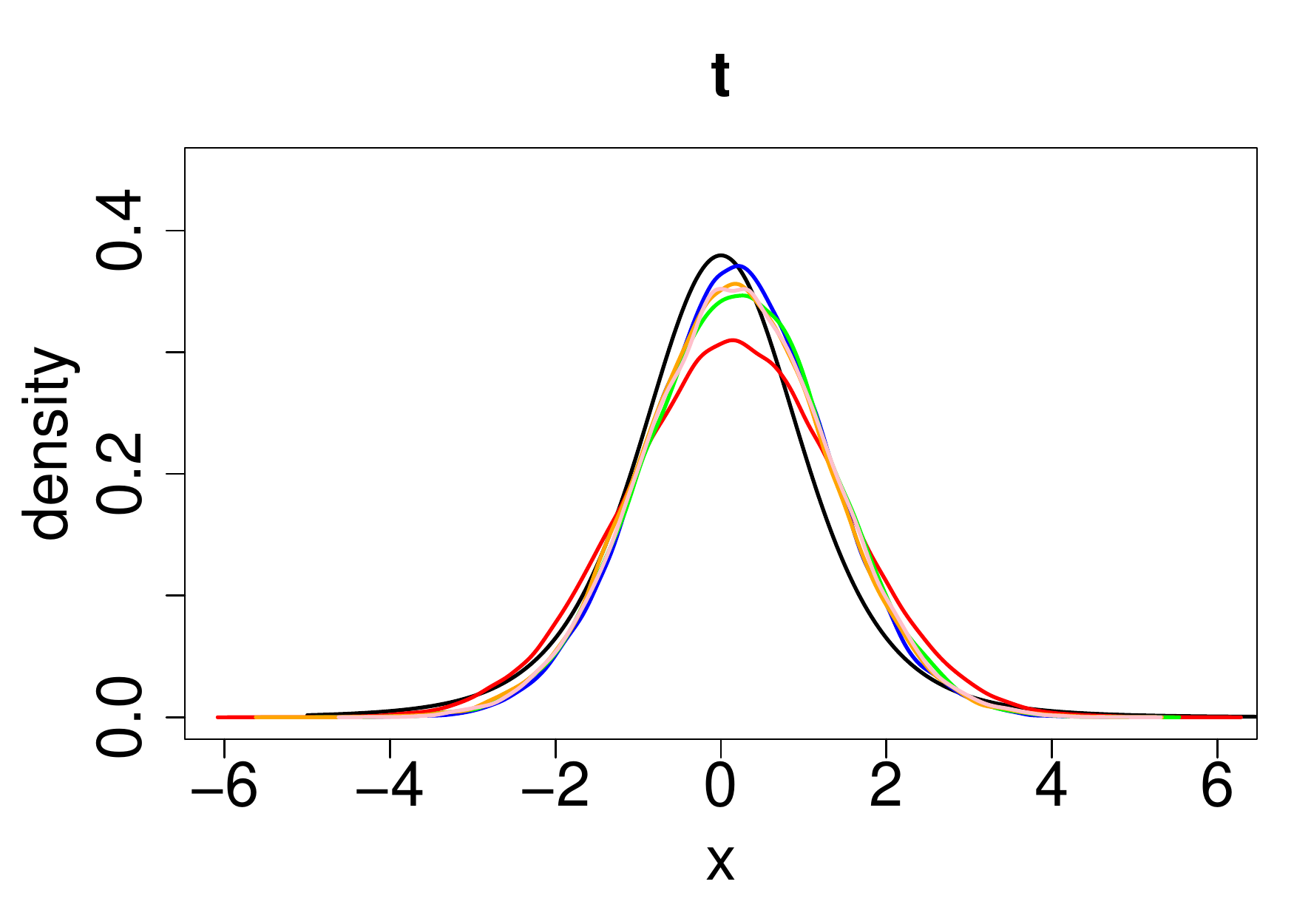}
\includegraphics[width=5cm]{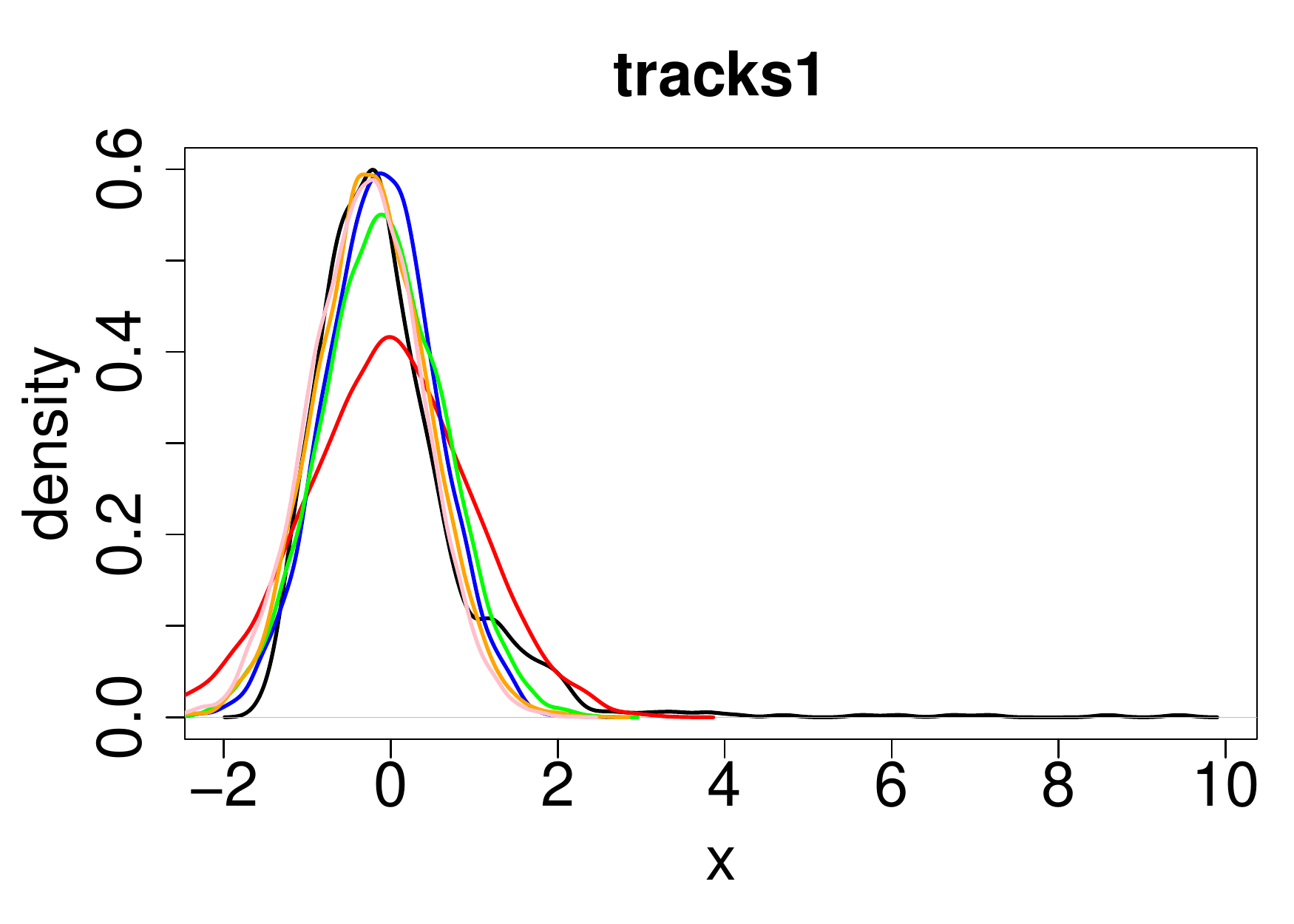}
\includegraphics[width=5cm]{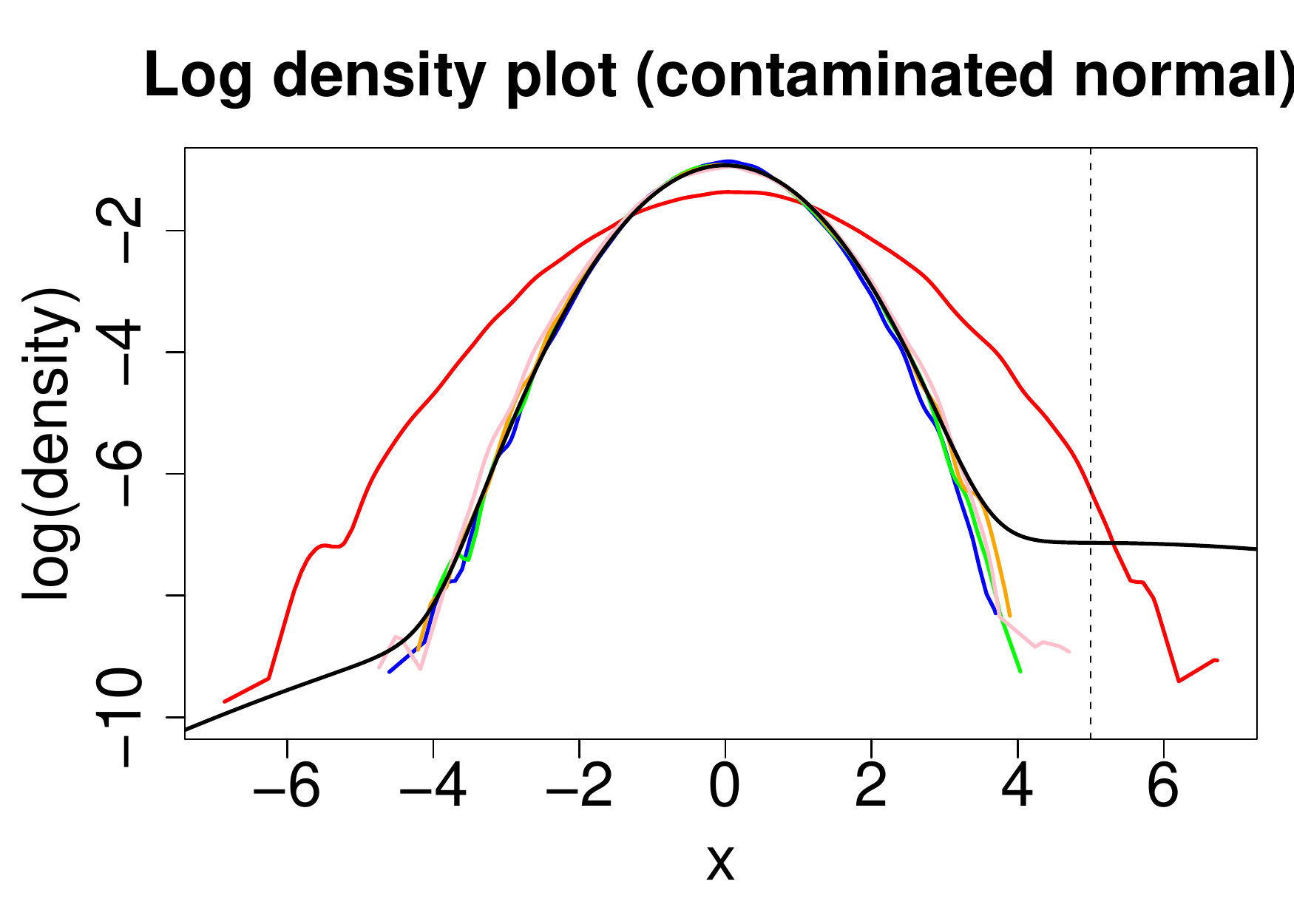}
\includegraphics[width=5cm]{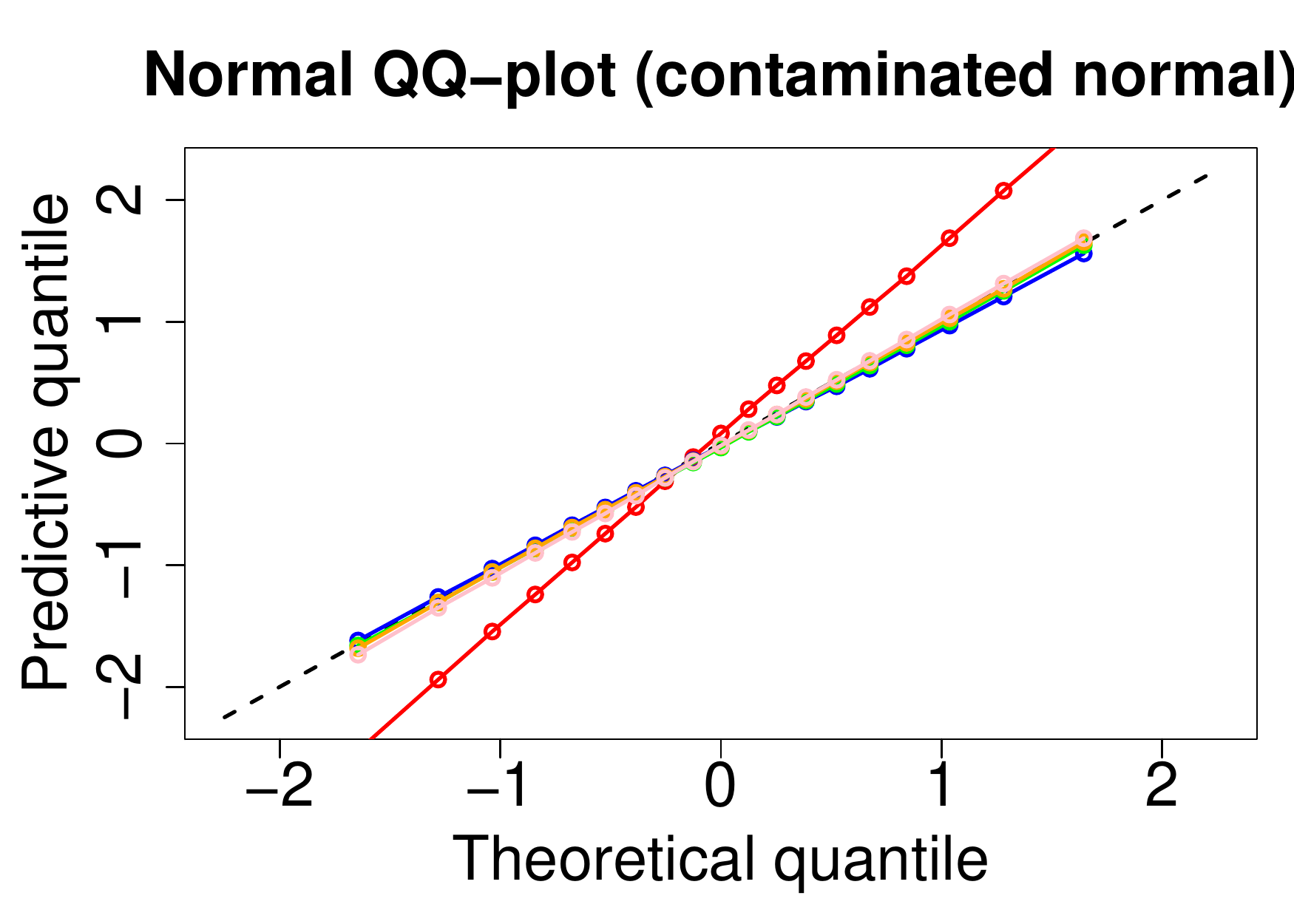}
\includegraphics[width=5cm]{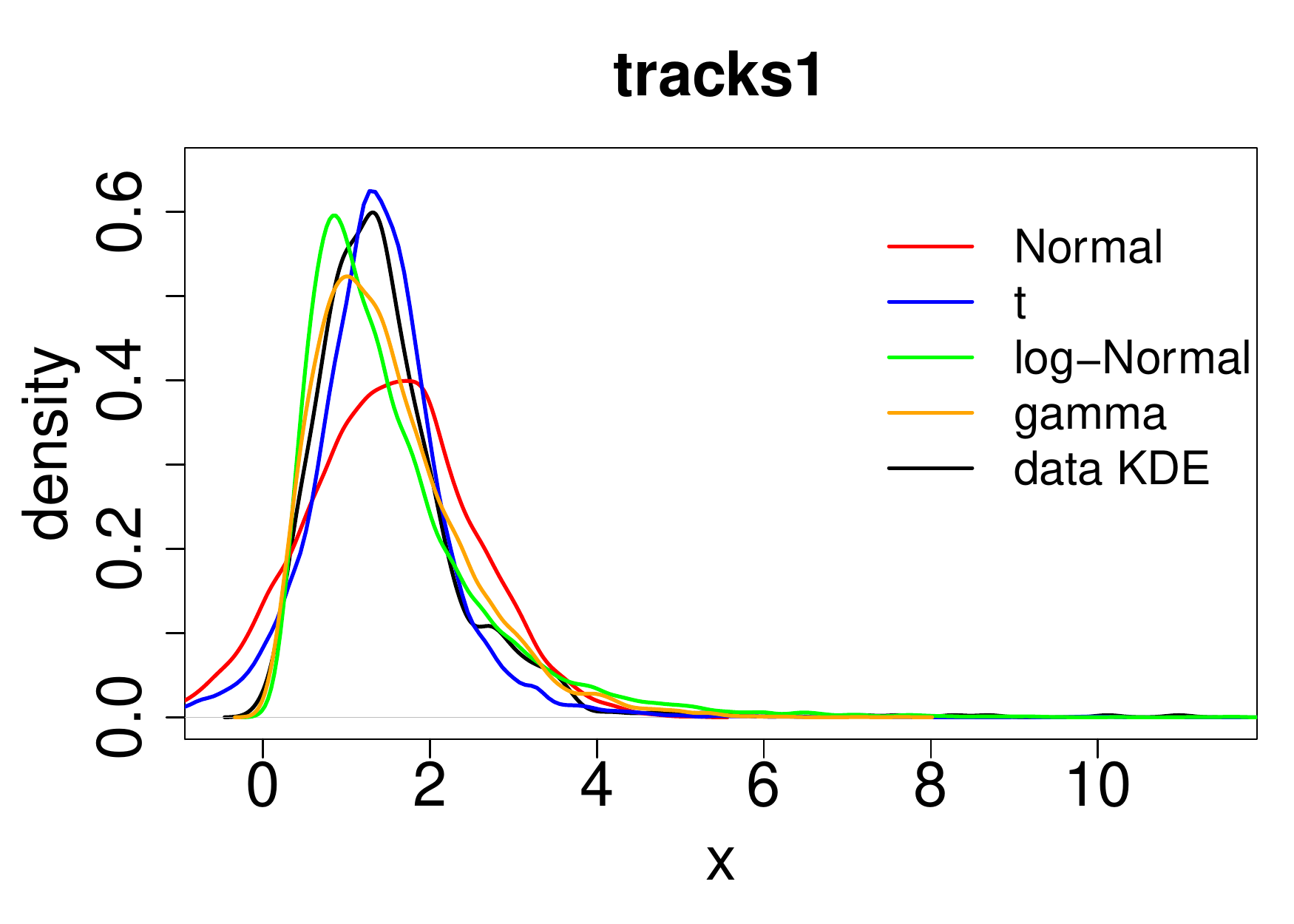}
\caption{\small Top: Posterior predictive distributions (smoothed from a sample) arising from Bayesian minimum divergence estimation fitting a normal distribution to an $\epsilon$-contaminated normal (left), a t-distribution (middle) and the tracks1 dataset (right) using the KL-Bayes (red), Hell-Bayes (blue), TV-Bayes (pink), alpha-Bayes (green) and power-Bayes (orange). Bottom: Log-density plots (smoothed from a sample) (left) and Normal QQ plots (middle) for the $\epsilon$-contaminated normal dataset. Lower right plots the posterior predictive distributions (smoothed from a sample) from alternative models using the KL-Bayes, Normal (red), t (blue), logNormal (green) and gamma (orange).}
\label{Fig:SynthDataGraphs}
\end{figure}

The experiments below demonstrate the robustness of the general Bayesian update targeting KL-divergence (red), Hellinger-divergence (blue), TV-divergence (pink), alpha-divergence (green) and power-divergence (orange) (in future these may be referred to as KL-Bayes, Hell-Bayes, TV-Bayes, alpha-Bayes and power-Bayes respectively). For illustrative purposes we have fixed $\alpha=0.75$ for the alpha-divergence and $\alpha=0.5$ for the power-divergence. Figure \ref{Fig:SynthDataGraphs} plots the posterior predictive originating from fitting a normal model $f(\cdot;\theta)=\mathcal{N}(\mu,\sigma^2)$ to two simulated data sets and one `real' data set. For the first data set $n=1000$ data points were simulated from a normal $\epsilon$-contamination distribution
\begin{equation}
g=0.99\times\mathcal{N}(0,1)+0.01\times\mathcal{N}(5,5^2),
\end{equation}

for the second $n=200$ were simulated from a Student's t-distribution with degrees of freedom 4 and the real data set, tracks1, was taken as the first variable from the `Geographical Original of Music Data Set'\footnote{downloaded from https://archive.ics.uci.edu/ml/datasets/Geographical+Original+of+Music} containing $n=1059$ data points, where a KDE of the data appeared to be approximately normally distributed. Prior distributions $\mu\sim\mathcal{N}(0,10^2)$ and $\sigma\sim\mathcal{G}(0.001,0.001)$ were used for all examples.

It is easy to see from the plots in figure \ref{Fig:SynthDataGraphs} that minimising the KL-divergence, requires that the density of the outlying contamination centred at $x=5$ to be captured correctly. This is at the expense of capturing the density of the remaining 99\% of the data, this is especially clear in the log-density plot. The TV-Bayes, alpha-Bayes and power-Bayes appear to correctly capture the distribution for 99\% of the data. The boundedness of the TV-divergence and the alpha-divergence means that the contamination is not allowed to unduly affect the analysis, while the power-divergence is able to down-weight the influence of the contaminated points as they are ascribed low predicted probability under the model for 99\% of the data. The Hell-Bayes appears to fit too small a variance even for 99\% of the data. This is due to the small sample efficiency problems which will be discussed in Appendix \ref{Sub:Efficiency}.

The Student's t-distribution has consistently heavier tails than the normal distribution. Thus the top right hand plot of figure \ref{Fig:SynthDataGraphs} more clearly demonstrates the importance placed on tail misspecification by each method. The KL-divergence fits the largest variance to correctly capture these tails for the most extreme tail observations, the alpha-divergence is able to fit a slight smaller variance because of its bounded nature. However the variance of the alpha-Bayes predictive is still larger than the other methods produce because of the greater convexity of the scoring function depicted in figure \ref{Fig:alphaKLH_scorecomp}. Lastly Hell-Bayes, TV-Bayes and power-Bayes place the least weight on tail misspecification and are therefore able to fit a smaller variance and produce a predictive more closely resembling the data generating process for the majority of the data.

 The importance given to tail misspecification is also evident in the `tracks' example. Here there is an ordering from TV, power, Hellinger, alpha and KL divergence on both bias towards the right tail and on the size of the predictive variance, in response to the slight positive skew of the KDE of the data. Here it is clear to see that the power, TV, Hellinger and power divergence produce much better fits of the majority of the data than the other methods do. Lastly, the bottom right plot demonstrates how several possible alternative models to the Gaussian perform on the `tracks1' data set\footnote{the data set was transformed by adding $\min(tracks1)+0.001$ to every value in order to make the data strictly positive so the gamma and log-Normal distributions could be applied}, when updating using the KL-divergence. This shows that a Gaussian distribution was actually the best fit for the bulk of the data, and the poor fit achieved is down to the importance placed on tail misspecification by the estimating procedure rather than the model selected.

\subsection{Regression under heteroscedasticity}

In addition to the simple inference examples above we consider how changing the divergence can affect inference in a regression example. From the previous examples we can see that when the tails of the model are misspecified the KL-divergence minimising predictive distribution inflates the variance of the fitted model to ensure no observations are predicted with low probability. Further, placing large weight on tail observations, which occur with low probability, creates large variance across repeat sampling. This is exactly why parameter estimates in linear regression under heteroscedasticity errors have a large variance. 

Bayesian minimum divergence inference places less weight on tail observations: it is thus able to produce inferences with a smaller predictive variance and a smaller variance across repeat sampling. While repeat sampling and estimation variance are not problems in the Bayesian paradigm these results do show that traditional Bayesian inference can be somewhat imprecise when the tails are misspecified which is clearly undesirable when conditioning on observed data.

In order to demonstrate this we simulated $n=200$ data points with $N=50$ repeats from the following heteroscedastic linear model

\begin{align}
y\sim \mathcal{N}\left( X\boldsymbol{\beta},\sigma(X_1)^2\right), \textrm{ where } \sigma(X_1)=\exp\left(\frac{2X_1}{3}\right).
\end{align}

\begin{figure}[H]
\centering
\includegraphics[width=7cm]{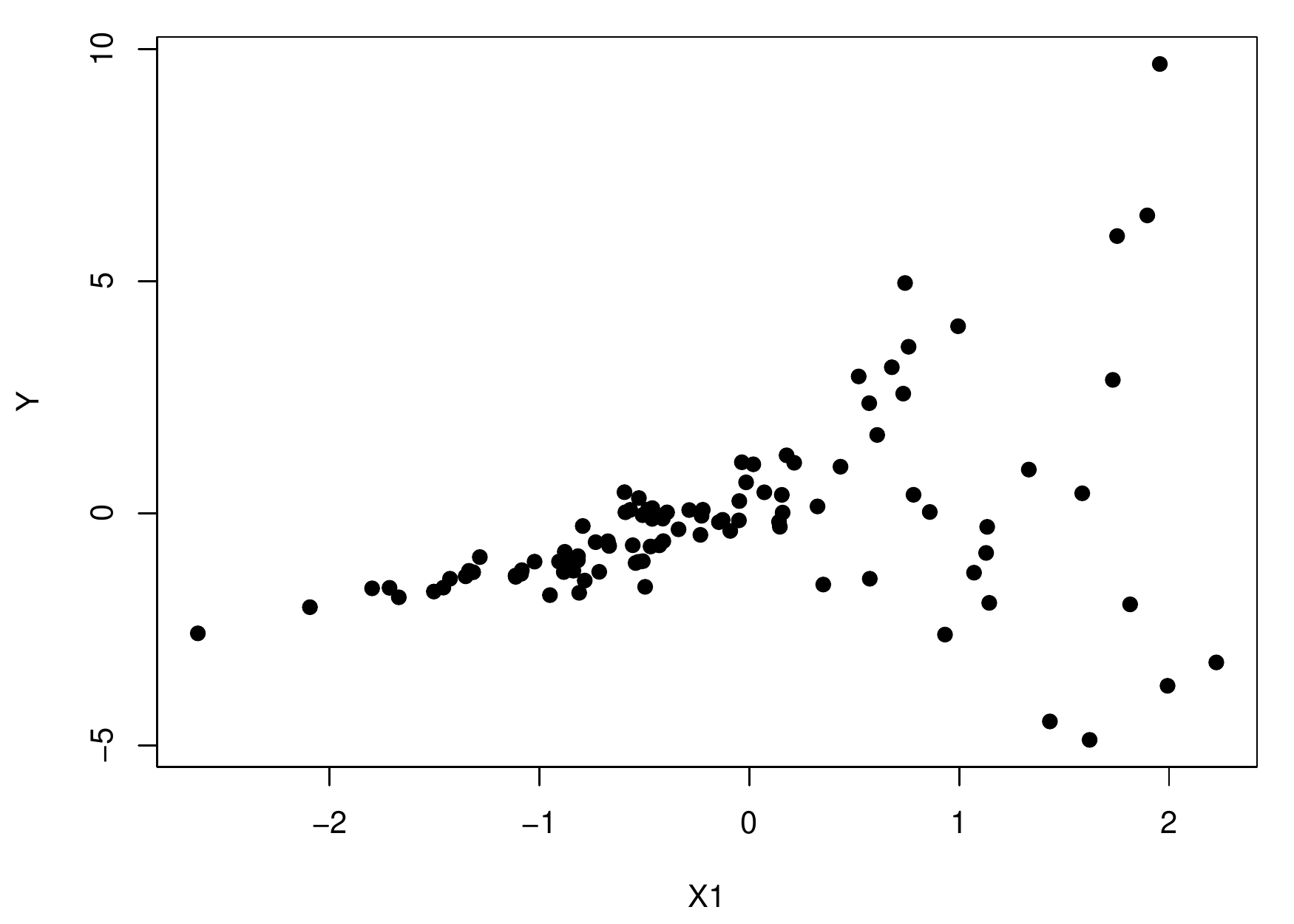}
\caption{\small One of the repeat data sets simulated from the heteroscedastic linear model with p=1.}
\label{Fig:heteroscedasticdata}
\end{figure}

For the experiments we simulated $X\sim\mathcal{N}_p\left(\mathbf{0},I\right)$ and each $\beta_i\sim\textrm{Unif}[-2,2]$ were held constant across experiments. Figure \ref{Fig:heteroscedasticdata} plots one realisation from this heteroscedastic linear model.

We then conducted general Bayesian updating using the five divergences mentioned above with priors $\sigma^2\sim\mathcal{IG}(2,0.5)$ and $\beta_i|\sigma^2\sim\mathcal{}(0,5\sigma^2)$. The mean squared errors (MSE) for the posterior means of the parameters $\frac{1}{N}\sum_{j=1}^n\sum_{i=1}^p\left(\hat{\beta}_i-\beta_i\right)^2$, the MSE for the predictive means on a test set of size 100 simulated from the model without error  $\frac{1}{N}\sum_{j=1}^n\sum_{i=1}^{100}\left(\hat{Y}_i-Y_i\right)^2$ and the predictive mean variances are presented in Table \ref{Tab:heteroscedastic} below. In order to apply the Hellinger, TV and alpha divergences to a regression problem an estimate of the conditional density of the response given the covariates is required. We follow authors \cite{hooker2014bayesian} and implement conditional KDEs to approximate the true data generating density. For simplicity, the familiar two-stage bandwidth estimation process of \cite{hansen2004nonparametric} was used to find the optimal bandwidth parameters.

 \begin{table}[ht]
 \small
 \centering
 \caption{\small Top: Table of posterior mean MSE values for parameters $\boldsymbol{\beta}$ and MSE for a test set simulated without error, across $N=50$ repeats from datasets of size $n= 200$ with the dimension of $\boldsymbol{\beta}$ $p=1,5,10,15,20$ under the Bayesian minimum divergence technology. Bottom: Table of posterior mean values for the variance of the linear model to be interpreted in terms of the precision of the predictive distribution minimising the respective divergence criteria.}
 \label{Tab:heteroscedastic}
\begin{tabular}{rrrrrrrrrrr}
  \hline\\[-0.9em]
     MSE & \multicolumn{2}{c}{\textbf{KL}} & \multicolumn{2}{c}{\textbf{Hell}} & \multicolumn{2}{c}{\textbf{TV}} & \multicolumn{2}{c}{\textbf{alpha}} & \multicolumn{2}{c}{\textbf{power}}\\
     \hline\\[-0.9em]
 & $\beta$ &  Test &  $\beta$ &  Test &  $\beta$ &  Test & $\beta$ & Test & $\beta$ & Test \\ 
\hline\\[-0.9em]
p=1 & 0.03 & 3.25 & 0.03 & 3.42 & 0.04 & 4.51 & 0.03 & 2.57 & 0.02 & 1.92 \\ 
   p=5 & 0.09 & 7.92 & 0.06 & 5.61 & 0.05 & 5.12 & 0.04 & 4.32 & 0.04 & 3.56 \\ 
   p=10 & 0.15 & 14.70 & 0.09 & 9.63 & 0.10 & 10.43 & 0.08 & 8.64 & 0.08 & 8.40 \\ 
   p=15 & 0.24 & 22.93 & 0.14 & 13.45 & 0.16 & 15.47 & 0.12 & 11.65 & 0.12 & 11.38 \\ 
   p=20 & 0.29 & 25.88 & 0.19 & 17.71 & 0.18 & 16.46 & 0.16 & 14.79 & 0.16 & 14.89 \\ 
    \hline
 \end{tabular}\\[0.5em]
 \begin{tabular}{rrrrrr}
   \hline\\[-0.9em]
  $\hat{\sigma}^2$ & \textbf{KL} & \textbf{Hell} & \textbf{TV} & \textbf{alpha} & \textbf{power} \\ 
   \hline\\[-0.9em]
 p=1 & 2.34 & 0.78 & 0.62 & 1.19 & 0.96 \\ 
   p=5 & 2.36 & 0.51 & 0.47 & 0.95 & 0.98 \\ 
   p=10 & 2.34 & 0.47 & 0.56 & 0.89 & 1.13 \\ 
   p=15 & 2.41 & 0.49 & 0.94 & 0.83 & 1.19 \\ 
   p=20 & 2.38 & 0.50 & 1.19 & 0.82 & 1.27 \\ 
   \hline
\end{tabular}
\end{table}

The bottom of table one demonstrates that the alternative divergences appear to be learning a smaller predictive variance than the KL-Bayes does under heteroscedastic errors. Under divergences alternative to the KL-divergence this is no longer an estimate of the variance of the responses in the data set, given the covariates. This should rather be interpreted predictively as the variance of the predictive distribution which is closest to the data generating distribution in terms of that alternative divergence. Therefore fitting a smaller variances demonstrates that the other divergences are placing more importance on fitting the majority of the data rather than just the outliers.

The top of table \ref{Tab:heteroscedastic} illustrates the impact fitting a large variance has on the parameter estimates of the mean function. Placing less influence on outliers allows all of the alternative divergences to produce more precise estimates of the parameters of the underlying linear relationship. This then leads to better performance when predicting the test set. Clearly the errors for all of the methods will increase as $p$ increases as the same amount of data is used to estimate more parameters. However it is clear that the errors under the KL-divergence are rising more rapidly. By being less sensitive to the error distribution the alternative divergences are better able to capture the true underlying process.

\subsection{Time series analysis}

In order to further demonstrate how inflating the variance by targeting the KL-divergence under misspecification can damage inference, we consider a less trivial time series example. We simulate $x_1,\ldots,x_T$ from an auto-regressive process of order L (AR(L)), and then consider additive independent generalised auto-regressive conditionally heteroscedastic of order (1,1) (GARCH(1,1)) errors, $e_1,\ldots,e_t$ with
\begin{align}
x_t&=\sum_{i=1}^L \mu_ix_{t-i}+\epsilon \textrm{ with } \epsilon\sim\mathcal{N}(0,\sigma^2)\\
e_t&=\psi_t\epsilon_k\textrm{ with } \epsilon\sim\mathcal{N}(0,1)\\
\psi_t^2&=\omega+\alpha_1e^2_{t-1}+\beta_1\psi^2_{t-1}\\
y_t&=x_t+e_t
\end{align}

where $\omega>0$, $\alpha_i>0$, $\beta_j\geq 0$. GARCH processes are used to model non-stationary, chaotic time series where the variance of the process depends on the magnitude and sign of the previous observations of the process. Eliciting a GARCH process from a DM is a difficult task. It is far from obvious how this GARCH process behaves as a function of its parameters and selecting a lag length for the AR process as well as two lag lengths for the GARCH Process increases the complexity of the model selection problem. Therefore it seems conceivable that a DM could want to fit a simple AR process to noisy time series data in order to investigate the underlying process. One situation where this may be desirable is in financial time series applications where large amounts of data can arrive at a very high frequency.

In order to investigate how the minimum divergence methods perform in this scenario we simulated 3 data sets with $T=1000$, and fitted and AR(L) process to these, where $L$ was chosen to match the underlying AR process. The 3 data sets were given by

\begin{enumerate}
\item an AR(3) with $\mu = (0.25,0.4,0.2,0.3)$
\item an AR(1) with $\mu=(0,0.9)$ with GARCH(1,1) errors $\omega = 2, \alpha_1 = 0.99, \beta_1 = 0.01$
\item an AR(1) with $\mu=(0,0.9)$ with GARCH(1,1) errors $\omega = 1, \alpha_1 = 0.75, \beta_1 = 0.01$
\end{enumerate}

\begin{figure}[H]
\centering
\includegraphics[width=5cm]{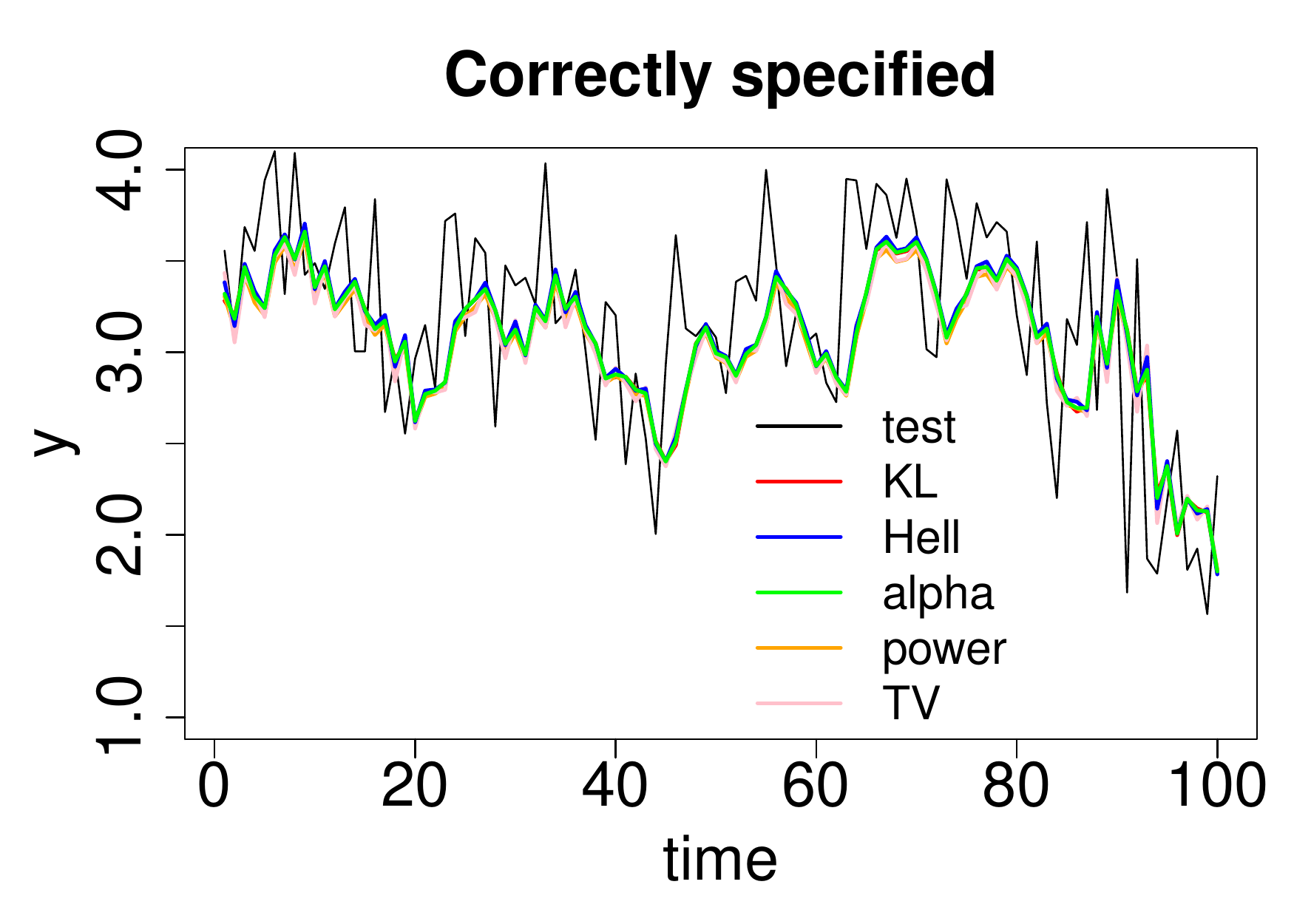}
\includegraphics[width=5cm]{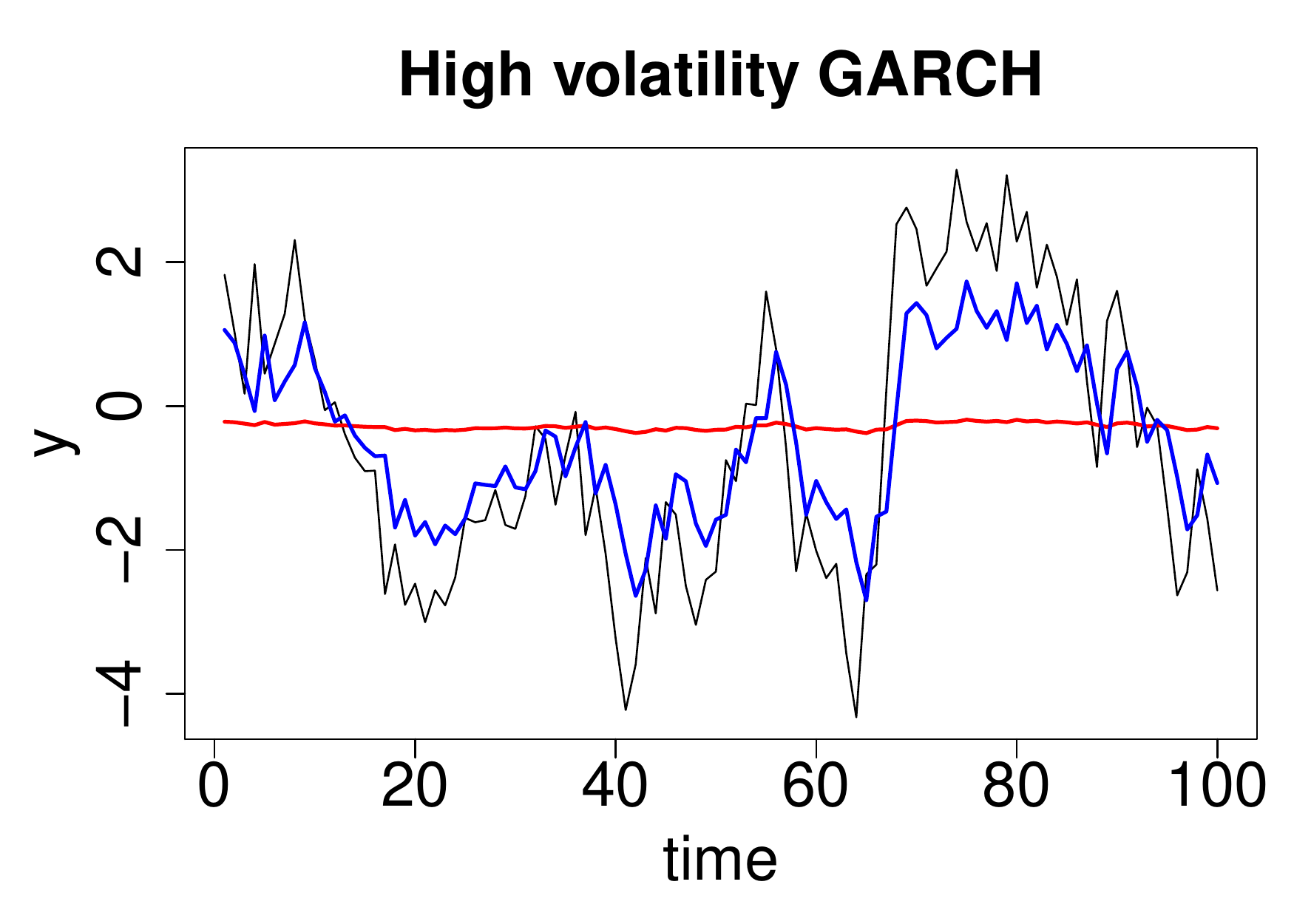}
\includegraphics[width=5cm]{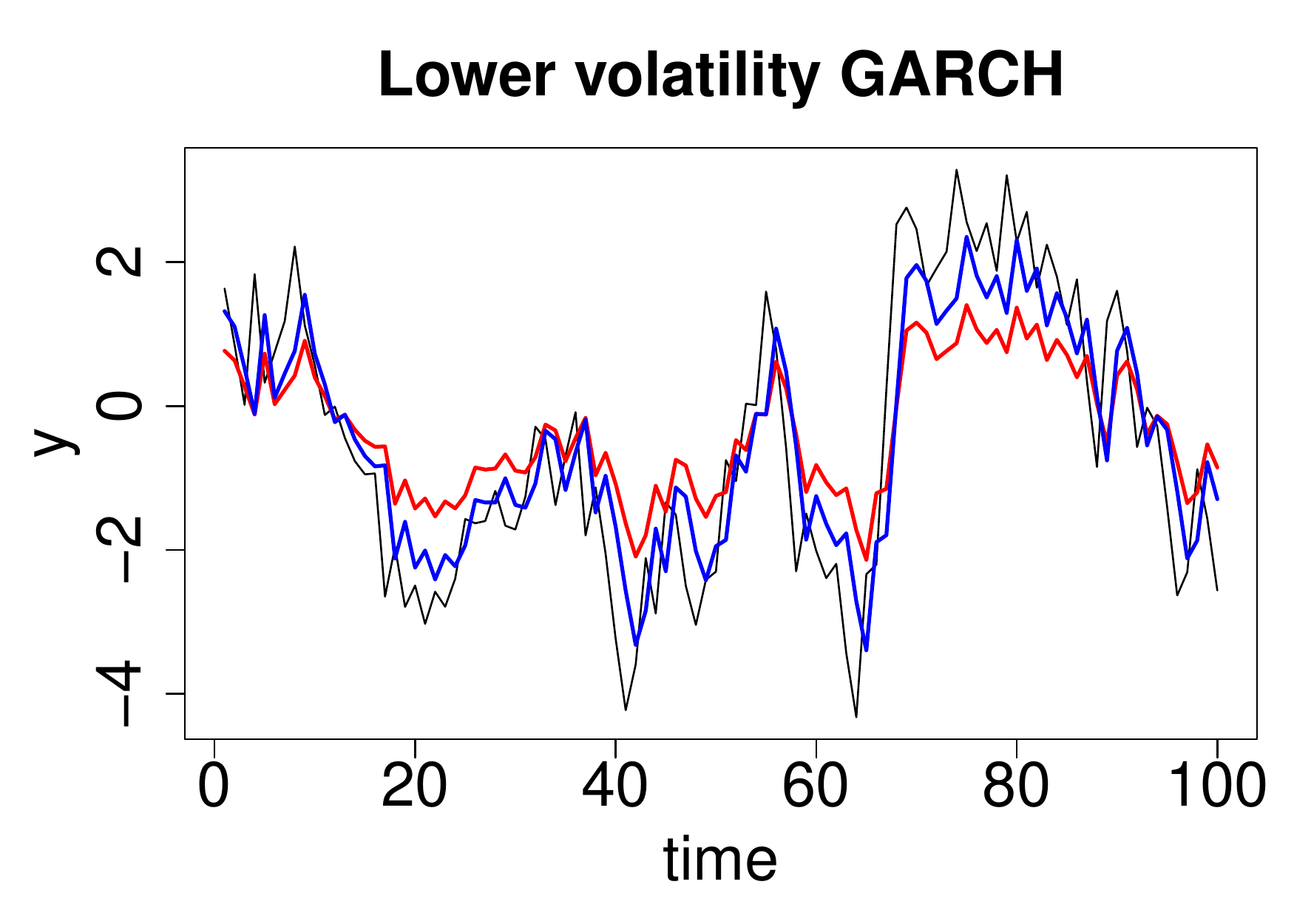}
\includegraphics[width=5cm]{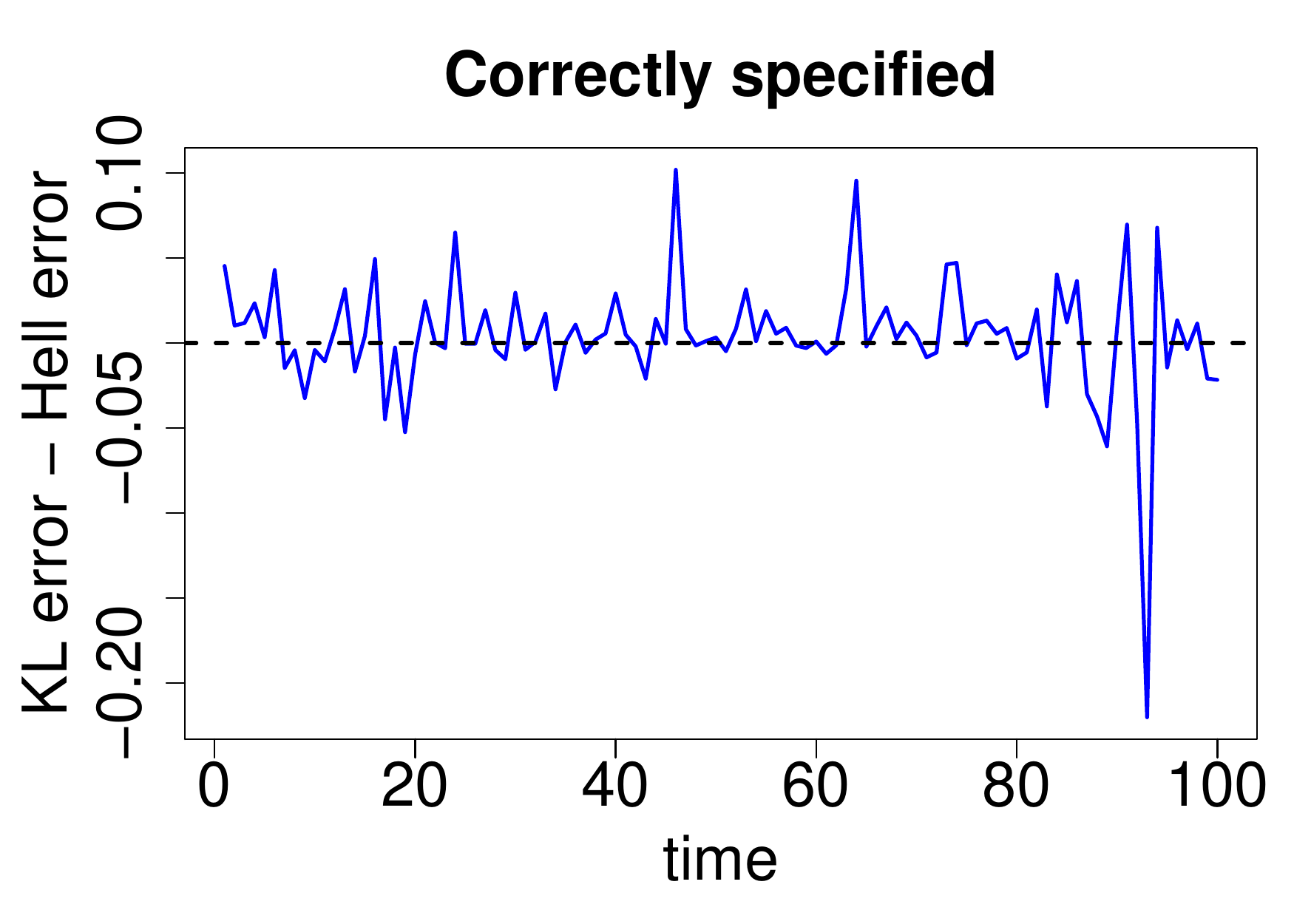}
\includegraphics[width=5cm]{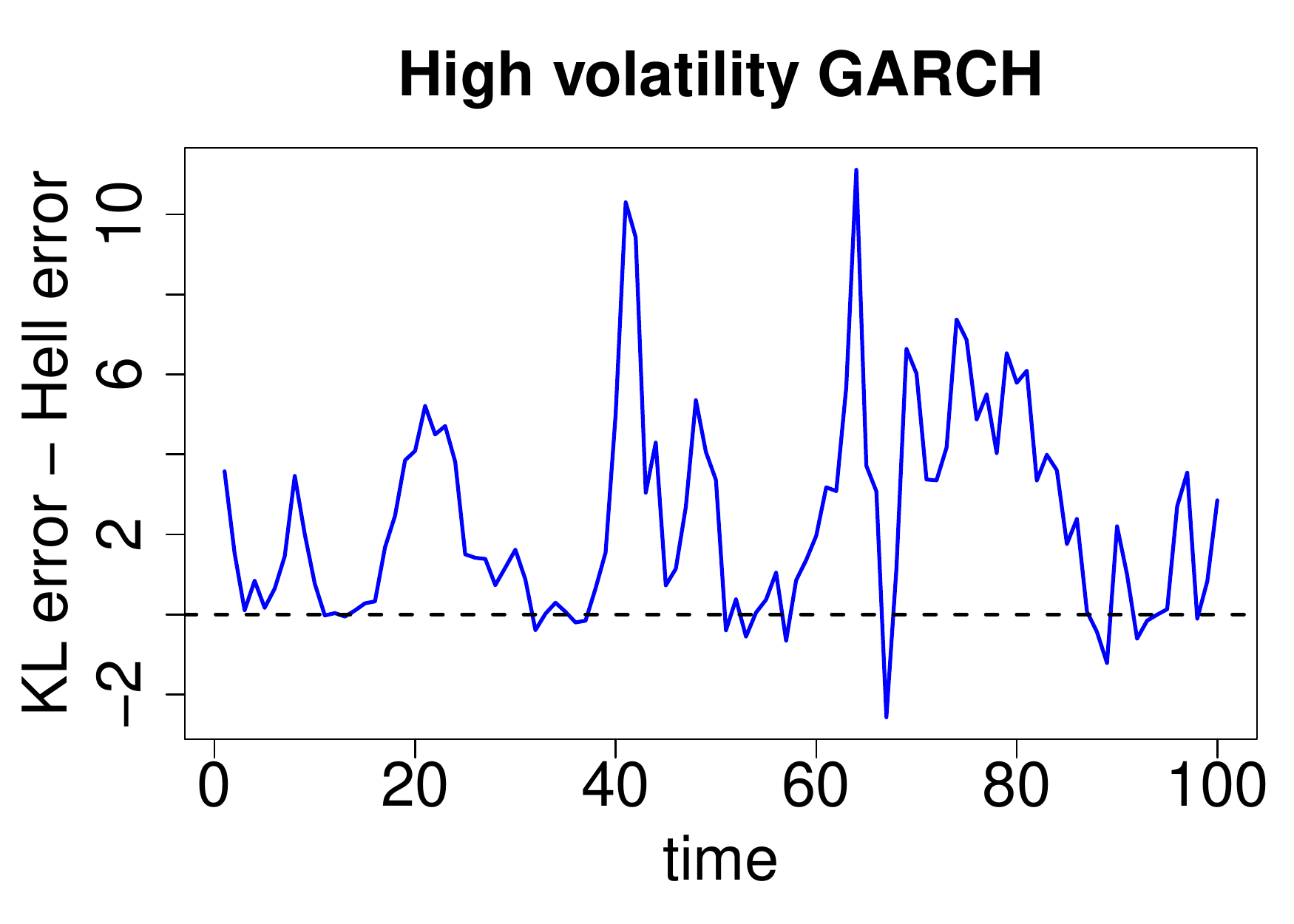}
\includegraphics[width=5cm]{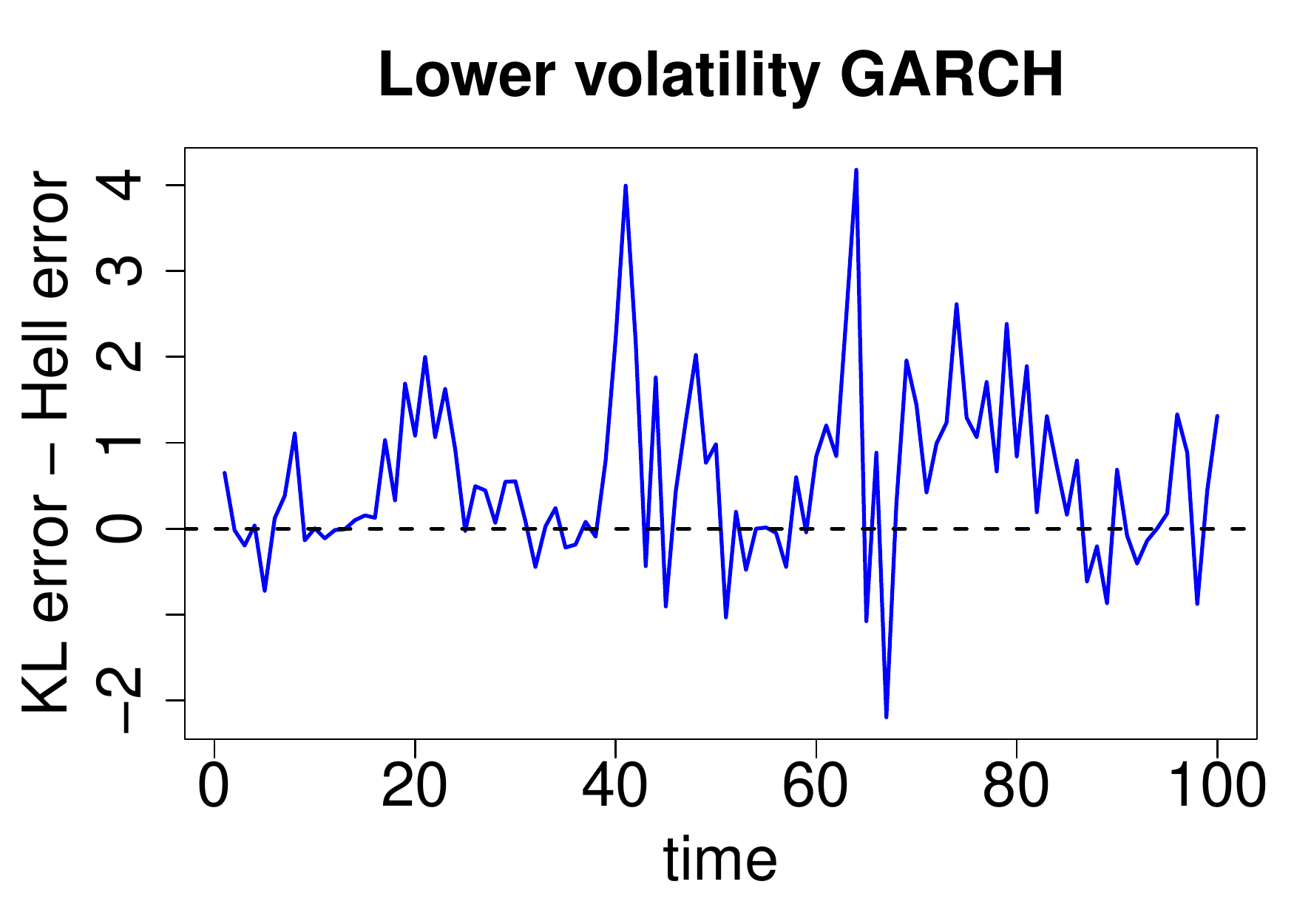}
\caption{\small Top: One step ahead posterior predictions arising from Bayesian minimum divergence estimation fitting to correct AR model to an AR(3) (left), an AR(1) with GARCH(1,1) errors, $\alpha=0.99$, $\omega =2$ (middle) and a AR(1) with GARCH(1,1) errors, $\alpha=0.75$, $\omega =1$ (right) using the KL-Bayes (red), Hell-Bayes (blue), TV-Bayes (pink), alpha-Bayes (green) and power-Bayes (orange). Bottom: the difference in one step ahead posterior squared prediction squared errors between the KL-Bayes and the Hell-Bayes.}
\label{Fig:TimeSeriesGraphs}
\end{figure}

The plots in figure \ref{Fig:TimeSeriesGraphs} demonstrate the one-step ahead posterior predictive performance of the minimum divergence posteriors on a test set $T=100$, simulated from the underlying AR process. Under misspecification we show only the inference under the Hell-Bayes to avoid cluttering the plots, the other minimum divergence posteriors perform similarly. We use the same priors as the regression example and once again conditional density estimates were used for the Hellinger, TV and alpha divergences.

The left hand plots demonstrate that when the model is correctly specified the minimum divergence posteriors produce similar time series inference to the KL-divergence. This is most easily seen in the lower plot which shows the difference in the squared prediction errors between the KL-Bayes and the Hell-Bayes is mostly around 0 and table \ref{Tab:RMSE} which demonstrates the root mean squared error across the test set is the same under both methods. The middle plots demonstrate how the KL-Bayes and the Hell-Bayes perform under `extreme' volatility in the error distribution. The tails of the model being very poorly specified causes the KL-Bayes posterior to fit a huge variance with the average posterior predictive variance across the test data set being slightly above 26. Fitting such a high variance makes the inference on the lag parameters $\mu$ insensitive to the data. As a result the underlying trend in the data is completely missed. In contrast the Hell-Bayes posterior predictive distributions have a much smaller variance of around 7. This allows the inference on the lag parameters to be much more sensitive to the underlying AR process. Clearly the Hell-Bayes is unable to exactly fit the truth as the model is misspecified. However it does a much better job of capturing the broad features of the underlying dependence between time points. This is clearly demonstrated by the considerably lower root mean squared predictive error showed in table \ref{Tab:RMSE} and by the error differences plot being mostly large and positive. The right hand plots demonstrates how the KL-Bayes and Hell-Bayes perform when the error distribution is less volatile. When the volatility is smaller the KL-Bayes predictive variance is also smaller. Therefore the inference is more sensitive to the underlying trend in the data than in the previous example. However the true dependence in the data is still under estimated relative to the Hell-Bayes, this is again demonstrated in table \ref{Tab:RMSE} and the error difference plot. Once again this shows that the way in which the KL-Bayes deals with misspecification, increasing the predictive variance, can mask some of the underlying trends in the data which can be discovered by other methods. Table \ref{Tab:RMSE} plots the root mean squared errors (RMSE) for the KL-Bayes and Hell-Bayes on the test set to quantify their correspondence to reality. 

\begin{table}[ht]
\small
 \caption{\small Root mean squared errors (RMSE) for the KL-Bayes and Hell-Bayes posterior mean predictions for 100 test data points from the underlying AR model.}
 \label{Tab:RMSE}
\centering
 \begin{tabular}{rrrr}
   \hline\\[-0.9em]
 RMSE & Correctly Specified & High Volatility & Low Volatility \\ 
   \hline\\[-0.9em]
 KL & 0.49 & 1.87 & 1.21 \\ 
   Hell & 0.48 & 1.07 & 0.94 \\ 
    \hline
 \end{tabular}
 \end{table}

\subsection{Application to high dimensions}\label{Sub:HighDim}

The examples is this paper only demonstrate the performance of these minimum divergence techniques for relatively small dimensional problems so as to clearly demonstrate the effect that tail misspecifications can have. However, it is as the dimension and complexity of the problem increases that these methods become more and more important, this can be put down to two aspects. The first of these is that outliers or highly influential contaminant data-points become hard to identify in high dimensions. In our examples it is clear from looking at the KDE or histogram of the data that there are going to be outlying observations, but in many dimensions visualising the data in this way is not possible. In addition to this automatic methods for outlier detection struggle in high dimensions \citep{filzmoser2008outlier}. 

The second factor in requiring robust inference in high dimensions, is that not only are outliers harder to identify, they are more likely to occur. The occurrence of outliers could be reinterpreted as indicating that the DM's belief model is misspecified in the tails. The fact that these occur should be unsurprising. We have already discussed specifying beliefs about tail behaviour requires thinking about very small probabilities which is known to be difficult to do \cite{winkler1968evaluation,o2006uncertain}, and often routine assumptions (for example Gaussianity) may be applied. As the dimension of the space increase the tails of the distribution account for a greater proportion of the overall density, increasing the chance of seeing observations that differ from the practitioners beliefs.

\section{Discussion}

This article uses general Bayesian updating \cite{bissiri2016general} in order to theoretically justify a Bayesian update that targets the parameters of a model that minimise a statistical divergence to the data generating process that is not the KL-divergence. When the $M$-open world is considered, moving away from targeting the minimisation of the KL-divergence can provide an important tool in order to robustify a statistical analysis. The desire for robustness ought to only increase as increasingly bigger models are built to approximate more complex real world processes. This paper provides the statistical practitioner with the principled justification to select the divergence they use for their analysis in a subjective manner allowing them the potential to make more useful predictions from their best approximate belief model. 

As it stands however more work is required to advise on the selection of the divergence used for the updating. Ways to more clearly articulate the impact of choosing a certain divergence to a DM needs further research along with ways to help guide the choice of any hyperparameters associated with the divergence. Further experimentation with complex real world data sets is also required to analyse how this robustness-efficiency trade-off associated with the selected divergence manifests itself in practice.

Another issue not addressed in this paper is how to tailor computational algorithms to the inferences we describe here. There exists a vast literature on optimising MCMC algorithms to sample from traditional Bayesian posteriors and in order to fully take advantage of the subjectivity this paper allows a statistician, a whole new class of computational algorithms tailored to different divergences may be required. In addition several of the divergences mentioned in this paper require a density estimate of the underlying process and further research into effectively doing this for complex high dimensional datasets can only improve the performances of these methods for real world problems.




\vspace{6pt} 




\subsubsection*{Acknowledgments}
We are indebted to 3 anonymous referees in helping us revise the structure and examples in this article. The authors would also like to thank Jeremias Knoblauch for his helpful discussions when revising the article. The first author acknowledges his research fellowship through the OxWaSP doctoral programme, funded by EPSRC.

\bibliography{bib}





\appendix
\section{M-closed efficiency}\label{Sub:Efficiency}


The examples in Section \ref{Sec:Illustrations}, 
demonstrate how minimising an alternative divergence to the KL-divergence can lead to superior robustness to tail misspecifications. We next illustrate that by placing less importance on tail misspecifications in order to gain improved robustness, the DM must trade-off a decrease in efficiency in the case when the model class does in fact contain the data generating process. Minimising the KL-divergence uses Bayes rule and thus conditions on the data coming from the model. As a result minimising the KL-divergence is bound to perform the best when the model class contains the data generating distribution. Here we illustrate this but also demonstrate the reassuring fact that the MDE methods perform comparatively even in this situation.

In order to examine the frequentist efficiency, the observed mean squared error (MSE) $\frac{1}{N}\sum_{j=1}^N(\hat{\theta}_n^j-\theta)^2$, over $N$ repeats of a simulated experiment are examined. The MSE's are examined on two examples $X_1\sim\mathcal{N}(0,10^2)$ and $X_2\sim\mathcal{N}(15,10^2)$ when fitting the model $X_i\sim\mathcal{N}(\mu_i,\sigma_i^2)$  with prior $\mu_i\sim\mathcal{N}(0,5^2)$ and $\sigma_i\sim\mathcal{G}(0.01,0.01)$ for $i=1,2$. The prior for $\mu_1$ is centred on the data while the prior for $\mu_2$ is not. Table \ref{Tab:MSEX} plots the observed MSEs over $N=50$ replications\footnote{The \textit{stan} output from a number of these experiments associated with the TV divergence suggested that the sample from the posterior was not sufficient to use for posterior inference, for the purpose of this preliminary analysis these experiments were disregarded here, but this suggests there exists computational issues associated with the TV-Bayes.}.

\begin{table}[H]
\centering
\small
\caption{\small Table of posterior mean MSE values for $\mu_i$, $\sigma_i$ from data sets of size $n= 50, 100, 200, 500$ under the Bayesian minimum divergence technology. }
\label{Tab:MSEX}
\begin{tabular}{rrrrrrrrrrrr}
   \hline\\[-1em]
  &MSE & \multicolumn{2}{c}{\textbf{KL}} & \multicolumn{2}{c}{\textbf{Hell}} & \multicolumn{2}{c}{\textbf{TV}} & \multicolumn{2}{c}{\textbf{alpha}} & \multicolumn{2}{c}{\textbf{power}}
\\   \hline\\[-1.3em]
   & & $\mu_1=0$ & $\mu_2=15$ & $\mu_1=0$ & $\mu_2=15$ & $\mu_1=0$ & $\mu_2=15$ & $\mu_1=0$ & $\mu_2=15$ & $\mu_1=0$ & $\mu_2=15$ \\
   \hline
 $\mu$ &n=50 & 1.89 & 4.43 & 1.50 & 16.88 & 2.68 & 6.77 & 2.03 & 4.29 & 0.76 & 93.80 \\ 
   &n=100 & 0.64 & 1.29 & 0.57 & 5.19 & 1.20 & 1.83 & 0.67 & 1.25 & 0.44 & 33.80 \\ 
   & n=200 & 0.35 & 0.47 & 0.32 & 1.50 & 0.48 & 0.50 & 0.35 & 0.47 & 0.31 & 7.83 \\ 
   & n=500 & 0.15 & 0.16 & 0.14 & 0.33 & 0.28 & 0.31 & 0.15 & 0.16 & 0.20 & 1.42 \\ \hline
   $\sigma$ & n=50 & 1.11 & 1.17 & 1.51 & 1.28 & 3.47 & 6.09 & 1.13 & 1.13 & 8.23 & 67.94 \\ 
   & n=100 & 0.60 & 0.61 & 0.98 & 0.87 & 1.30 & 1.61 & 0.65 & 0.64 & 1.99 & 10.64 \\ 
   & n=200 & 0.23 & 0.23 & 0.57 & 0.52 & 0.42 & 0.44 & 0.30 & 0.29 & 0.46 & 1.09 \\ 
  & n=500 & 0.10 & 0.10 & 0.26 & 0.26 & 0.29 & 0.27 & 0.12 & 0.13 & 0.18 & 0.21 \\ 
    \hline
 \end{tabular}
 \end{table}
 
Table \ref{Tab:MSEX} demonstrates several interesting points about these different methods. Firstly it appears as though the KL-Bayes, Hell-Bayes and alpha-Bayes perform very similarly when $n=500$, while the power-Bayes and TV-Bayes performs marginally worse. As we should expect, for $n<500$ the KL-Bayes appears to be optimal for both priors with the $\alpha$-Bayes appearing to perform the next best. Table \ref{Tab:MSEX} also shows the differing impacts the prior has on the different Bayesian updates. The KL-Bayes is least affected by the prior with the $\alpha$-Bayes and the TV-Bayes the next least. Authors \cite{hooker2014bayesian} recognised that the boundedness of the Hell-Bayes score function compared with the KL-Bayes score function, may cause problems such as this. The $\alpha$-Bayes score function is larger than the Hell-Bayes score function so our findings would be consistent with this. The power-Bayes score function with $\alpha=0.5$ is smaller than the Hell-Bayes score function when $g<0.25$. This explains why the power-Bayes is impacted more by the prior, and provides a likely explanation of why the power-Bayes appears to perform so well when the prior is specified around the data. For small values of $g$, which occur when the variance of $g$ is 100, the absolute size of the TV-Bayes score function is generally larger than the Hell-Bayes score, again explaining why the TV-Bayes performs better than the Hell-Bayes under a poorly specified prior.

Examining whether the estimates of $\sigma$ were generally above or below the data generating parameter demonstrates how these different methods learn, see table \ref{Tab:SignSEX}. For all sample sizes the KL-Bayes appears to over or under estimate the true variance in equal proportions over the 50 repeats, while the Hell-Bayes and alpha-Bayes, to a lesser extent, generally under predict the true variance. 

\begin{table}[H]
\small
 \centering
  \caption{\small Table of sums of posterior mean over (positive) or under (negative) estimation for $\mu_i$, $\sigma_i$ from data sets of size $n= 50, 100, 200, 500$ under the Bayesian minimum divergence technology}
\label{Tab:SignSEX}
 \begin{tabular}{rrrrrrrrrrrr}
   \hline\\[-1em]
  &MSE & \multicolumn{2}{c}{\textbf{KL}} & \multicolumn{2}{c}{\textbf{Hell}} & \multicolumn{2}{c}{\textbf{TV}} & \multicolumn{2}{c}{\textbf{alpha}} & \multicolumn{2}{c}{\textbf{power}} \\ 
   \hline\\[-1.3em]
   & & $\mu_1=0$ & $\mu_2=15$ & $\mu_1=0$ & $\mu_2=15$ & $\mu_1=0$ & $\mu_2=15$ & $\mu_1=0$ & $\mu_2=15$ & $\mu_1=0$ & $\mu_2=15$ \\
   \hline
 $\mu$ &n=50 & -12 & -32 & -10 & -50 & -5 & -37 & -12 & -32 & -8 & -50 \\ 
    &n=100 & -16 & -40 & -14 & -50 & -15 & -22 & -18 & -40 & -12 & -50 \\ 
    &n=200 & -2 & -24 & -8 & -46 & 1 & -10 & -2 & -26 & -2 & -50 \\ 
    &n=500 & -6 & -16 & -6 & -42 & 1 & 3 & -10 & -18 & -8 & -48 \\ \hline
   $\sigma$ &n=50 & 4 & 6 & -26 & -10 & 25 & 27 & -8 & -6 & 46 & 50 \\ 
    &n=100 & 0 & 0 & -36 & -32 & 19 & 12 & -14 & -14 & 38 & 48 \\ 
    &n=200 & -2 & -2 & -40 & -40 & 23 & 22 & -24 & -24 & 28 & 40 \\ 
    &n=500 & -2 & -2 & -36 & -40 & 27 & 21 & -28 & -28 & 14 & 20 \\ 
    \hline
 \end{tabular}
 \end{table}

\begin{figure}[H]
\centering
\includegraphics[width=3.8cm]{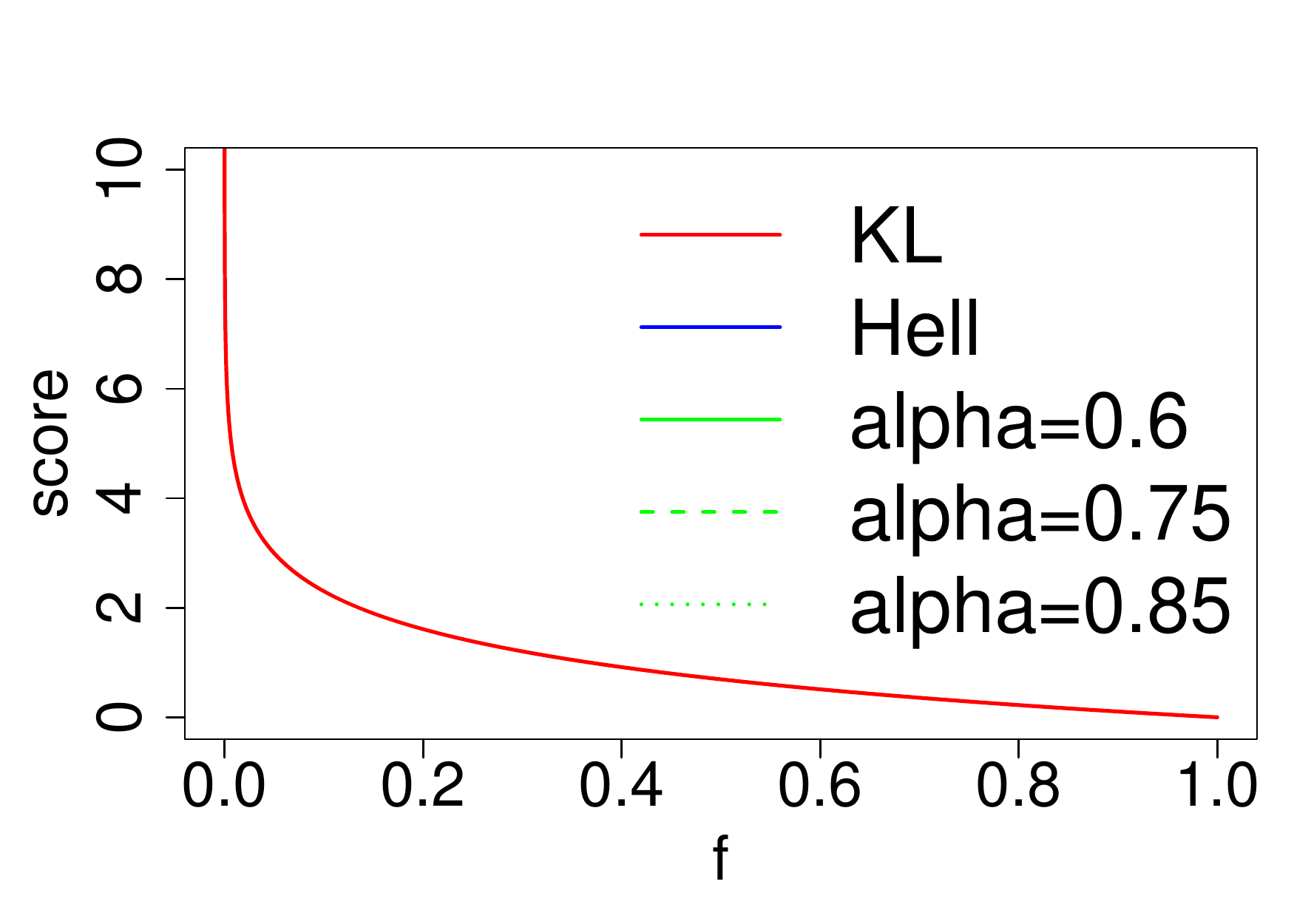}
\includegraphics[width=3.8cm]{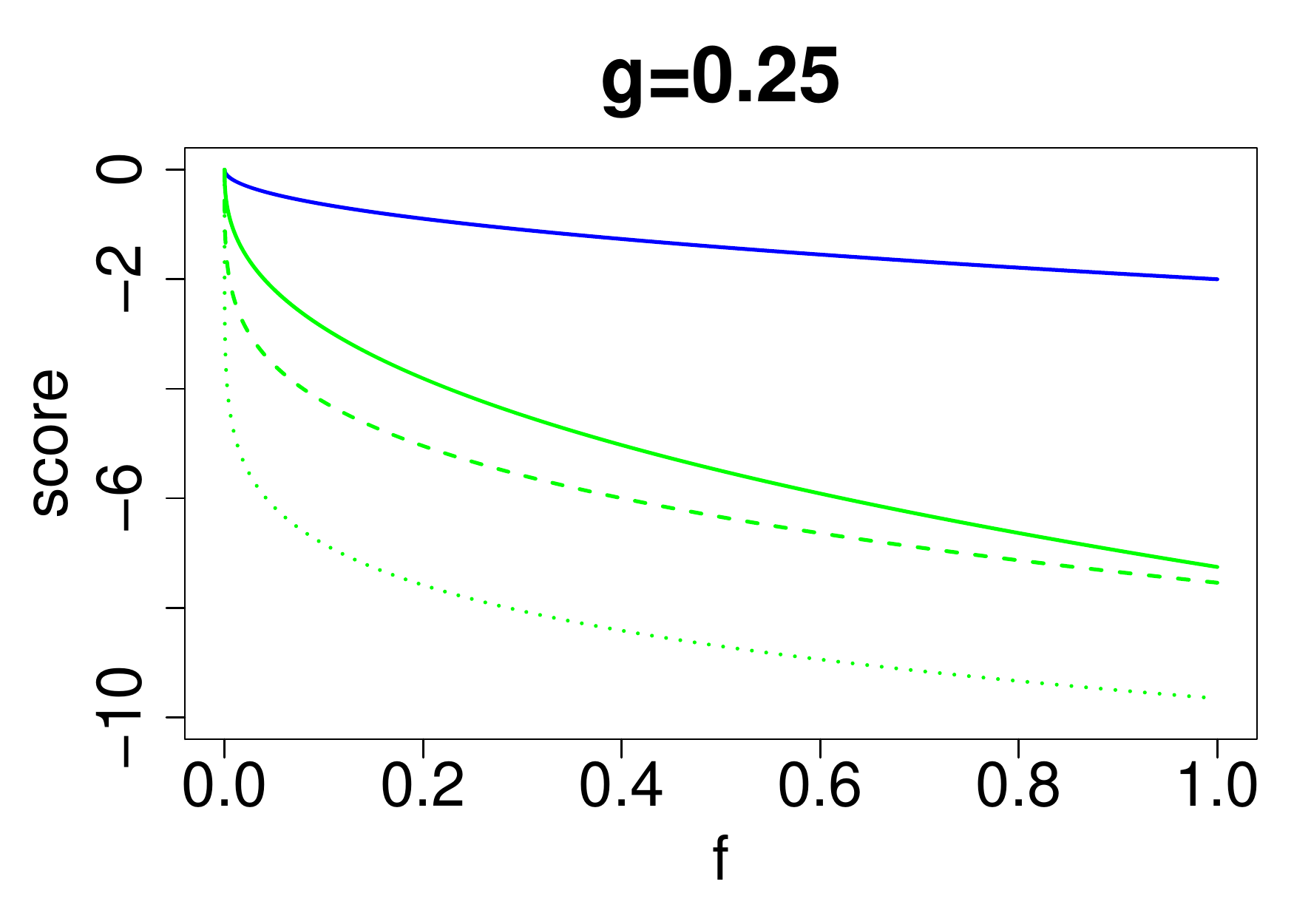}
\includegraphics[width=3.8cm]{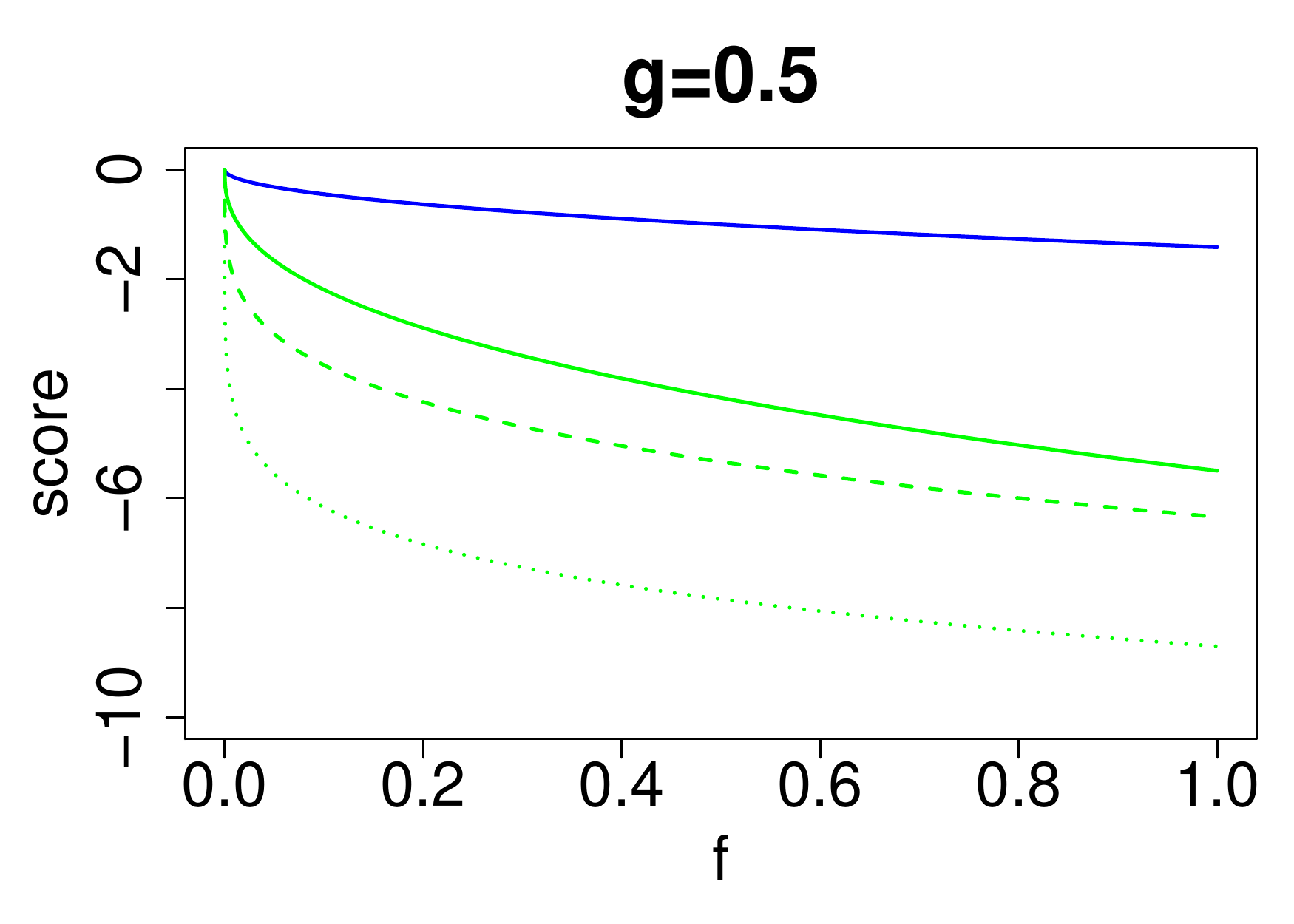}
\includegraphics[width=3.8cm]{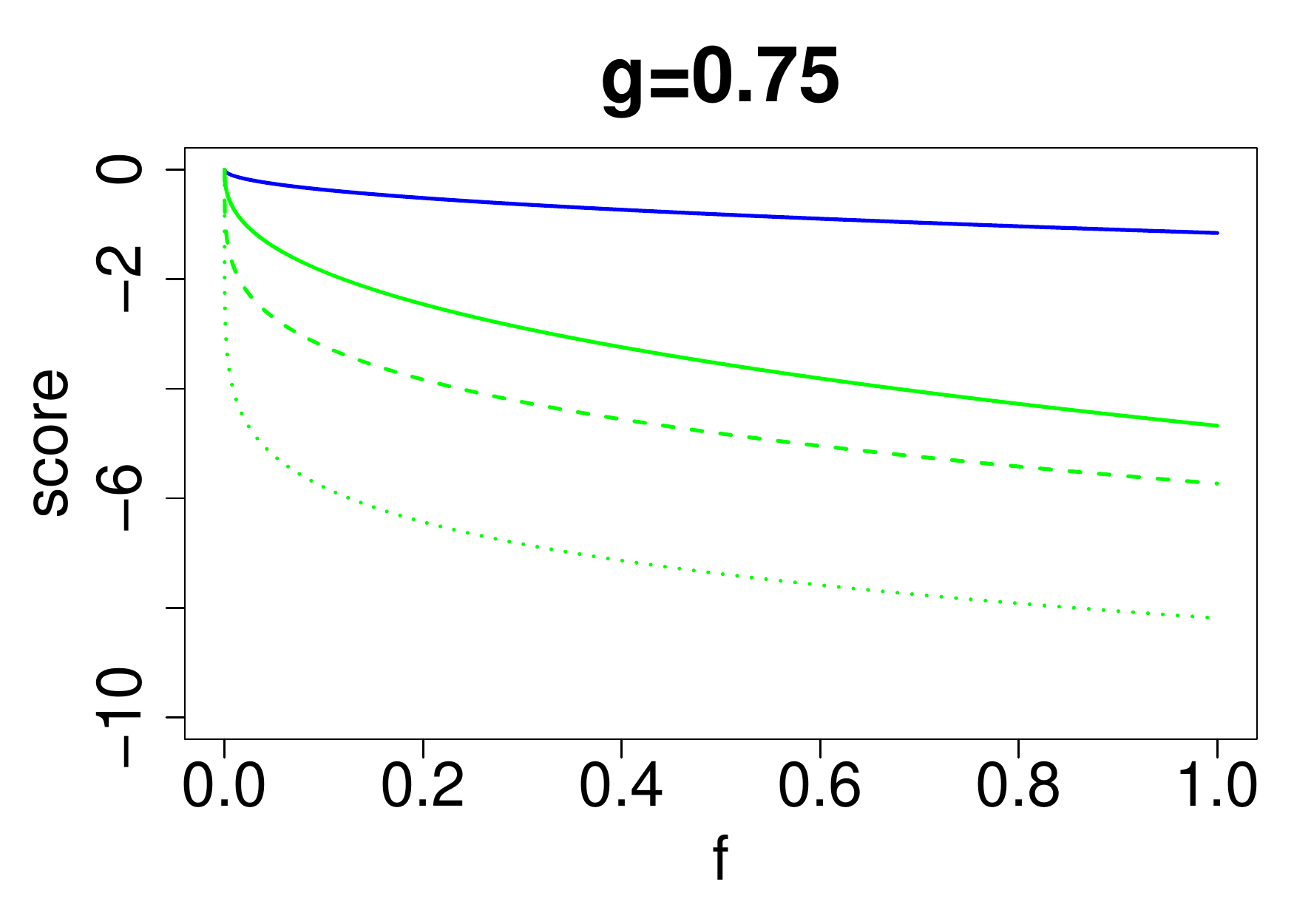}
\caption{\small Score comparisons when quoting probability $x$ for event with probability $g$. Red = KL, blue = Hellinger ($\alpha=0.5$, then divide by $4$), green = alpha ($\alpha=0.6$ = solid, $\alpha=0.75$ = dashed, $\alpha=0.85$ = dotted)}
\label{Fig:alphaKLH_scorecomp}
\end{figure}

Figure \ref{Fig:alphaKLH_scorecomp} compares the scoring functions used in the general Bayesian update targeting the KL, Hellinger and alpha divergences. As $\alpha$ increases, the upper bound on the score associated with the alpha-divergence increases and the scoring functions become more convex, tending towards the logarithmic score when $\alpha=1$. In contrast, the score associated with the Hellinger-divergence is closer to being linear. We believe that this is responsible for causing the smaller variance estimates when minimising the Hellinger-divergence and alpha-divergence when $\alpha=0.75$. The closer the score is to being linear the more similar the penalty is for over and under predicting the probability of an observation when compared with reality. The bounded nature of the scores prevents the penalty for under predicting being too high, and therefore more posterior predictive mass is placed near the MAP of the posterior predictive distribution. In contrast, a score function with greater convexity will penalise under prediction compared with the data generating density to a greater extent spreading the posterior mass out further. However, one must keep in mind that it is the sever nature of the penalty for under prediction incurred by the KL-divergence that renders it non-robust. On the other hand the TV-Bayes and the power-Bayes generally over predict the variance. For the power-Bayes this appears to be because the second term of the power divergence score is

\begin{equation}
\frac{1}{1+\alpha}\int f(y;\theta)^{\alpha+1}dy=\frac{1}{(1+\alpha)^{1+d/2}(2\pi)^{d\alpha/2}|\Sigma|^{\alpha/2}}
\label{Equ:powerDivRegTerm}
\end{equation}

for a multivariate Gaussian model, and this is decreased when the variances are made large. As was established in Section \ref{Sub:alphabetaDivergence}, the power-Bayes with $\alpha=0.5$ decreases the influence observations with small likelihood under the model have on the inference. In the example above the data has a variance of $\sigma^2=100$ and therefore the likelihood of each observation is very small. As a result the data is then excessively down-weighted by the power-Bayes when $\alpha=0.5$. When the sample size is small, a small increase in the variance may only decrease the observed likelihood term of the power-divergence score by a fraction and may increase the term in equation (\ref{Equ:powerDivRegTerm}) by more, resulting in a variance estimate that is too large. The score function for the TV-Bayes is the only score function considered here that is not monotonic in the predicted probability of each observation. While under the other scores, the score for each individual observation can be increased by predicting that observation with greater probability, this is not the case for the TV-Bayes. We believe this causes the TV-Bayes to produces variance estimates that are too large. We note that all of the scores remain proper for large samples and that these are just observations about how the different methods perform for small samples sizes when the model is correct.

\end{document}